%% file: MainPaper.tex
\documentclass[review,onefignum,onetabnum,dvipsnames]{siamart171218}

\input{SharedInfo}

\input{macros}

\usepackage{array,multirow,booktabs} % for fancier tables
\usepackage[export]{adjustbox} % To crop figures

\begin{document}
\nolinenumbers
\maketitle
\input{Abstract}

\input{KeywordsAms}

\input{Introduction}

\input{Method}

\input{Algorithm}

\input{Results}

\input{Conclusions}

\input{Acknowledgments}

\bibliographystyle{siamplain}
\bibliography{references}

\end{document}

%% file: SharedInfo.tex
\usepackage{amsfonts, amssymb,bbm }
\usepackage{graphicx}
\usepackage{mathrsfs,mathtools}
\usepackage{epstopdf}
\usepackage{algorithm}
\usepackage{algorithmic}
\usepackage{cleveref}
\usepackage{enumitem}

\setlist[enumerate,1]{label=(\arabic*),ref=(\arabic*)}

\ifpdf
  \DeclareGraphicsExtensions{.eps,.pdf,.png,.jpg}
\else
  \DeclareGraphicsExtensions{.eps}
\fi

\newcommand{\mathd}{\mathrm{d}}
\newcommand{\R}{\mathbbm{R}}
\DeclareMathOperator*{\argmin}{argmin}

\theorembodyfont{\upshape}

\newsiamremark{remark}{Remark}
\newsiamremark{hypothesis}{Hypothesis}
\crefname{hypothesis}{Hypothesis}{Hypotheses}
\newsiamthm{claim}{Claim}

\headers{Structure-preserving Nonlinear Filtering: CG \& DG}{Vidhi Zala, Robert M. Kirby, and Akil Narayan}

\title{Structure-preserving Nonlinear Filtering for Continuous and Discontinuous Galerkin Spectral/hp Element Methods\thanks{Submitted to the editors DATE.
}}

\author{Vidhi Zala\thanks{Scientific Computing and Imaging Institute and School of Computing, University of Utah, Salt Lake City, UT 84112
  (\email{vidhi.zala@utah.edu}\url{}).}
  \and Robert M. Kirby\thanks{
Scientific Computing and Imaging Institute and School of Computing, University of Utah, Salt Lake City, UT 84112
  (\email{kirby@cs.utah.edu}\url{ }).}
\and Akil Narayan\thanks{
Scientific Computing and Imaging Institute and Department of Mathematics, University of
Utah, Salt Lake City, UT 84112
  (\email{akil@sci.utah.edu}\url{}).}
}

\usepackage{amsopn}

%% file: macros.tex
\newcommand{\dx}[1]{\mathd #1}

\newcommand{\bs}[1]{\boldsymbol{#1}}
\renewcommand{\hat}[1]{\widehat{#1}}

\theoremstyle{example}

%% file: Abstract.tex
\begin{abstract}

Finite element simulations have been used to solve a variety of partial differential equations (PDEs) that model physical, chemical, and biological phenomena. The resulting discretized solutions to PDEs often do not satisfy requisite physical properties, such as positivity or monotonicity. Such invalid solutions pose both modeling challenges, since the physical interpretation of simulation results is not possible, and computational challenges, since such properties may be required to advance the scheme. We, therefore, consider the problem of computing solutions that preserve these structural solution properties, which we enforce as additional constraints on the solution. We consider in particular the class of convex constraints, which includes positivity and monotonicity. By embedding such constraints as a postprocessing convex optimization procedure, we are able to compute solutions that satisfy general types of convex constraints. For certain types of constraints (including positivity and monotonicity), the optimization is a filter, i.e., a norm-decreasing operation. We provide a variety of tests on one-dimensional time-dependent PDEs that demonstrate the efficacy of the method, and we empirically show that rates of convergence are unaffected by the inclusion of the constraints.

\end{abstract}

%% file: KeywordsAms.tex
% REQUIRED

% REQUIRED

%% file: Introduction.tex
\section{Introduction}\label{sec:intro}

Since the advent of numerical computing methods such as the finite element method (FEM) and the finite volume method (FVM) that solve partial differential equations (PDEs),  the scientific computing community has advanced these methods with the goal of having computed solutions that emulate real-world phenomena. Many such numerical methods rely on piecewise polynomial approximations of fields.  For example, we frequently numerically solve a system of PDEs to predict a state variable $u$ that advects with wave-like motion. When computationally representing $u$ using a piecewise polynomial, some physical properties of this variable can be lost sometimes. If $u$ denotes the concentration of some quantity, piecewise polynomial representations often  result in non-negative values $u$, yet negative values of $u$ are not physically interpretable.  Other examples of properties that may be lost are minimum/maximum bound on values or monotonically increasing or decreasing behavior. In this paper, we consider a general class of convex \textit{structural} properties of the state, which includes the previous examples. 

PDE models in which the importance of structure preservation can be observed include combustion problems, fluid flow problems, atmospheric predictions, and modeling of biochemical processes, such as platelet aggregation and blood coagulation. 

In many cases, the theoretical and practical feasibility of numerical methods depends on how closely the computed approximation to $u$ follows the requisite physical structure. If this underlying structure is violated, the resulting computation may produce unphysical predictions, and/or may cause solvability issues in numerical schemes that approximate PDE solutions. A violation in structure may arise from seemingly benign approximation properties; for example, polynomial approximations that yield Gibbs' oscillations often still converge in the mean-square sense, but these oscillations can cause a violation of positivity.

\begin{figure}[htbp] \centering
    \includegraphics[width=0.49\textwidth]{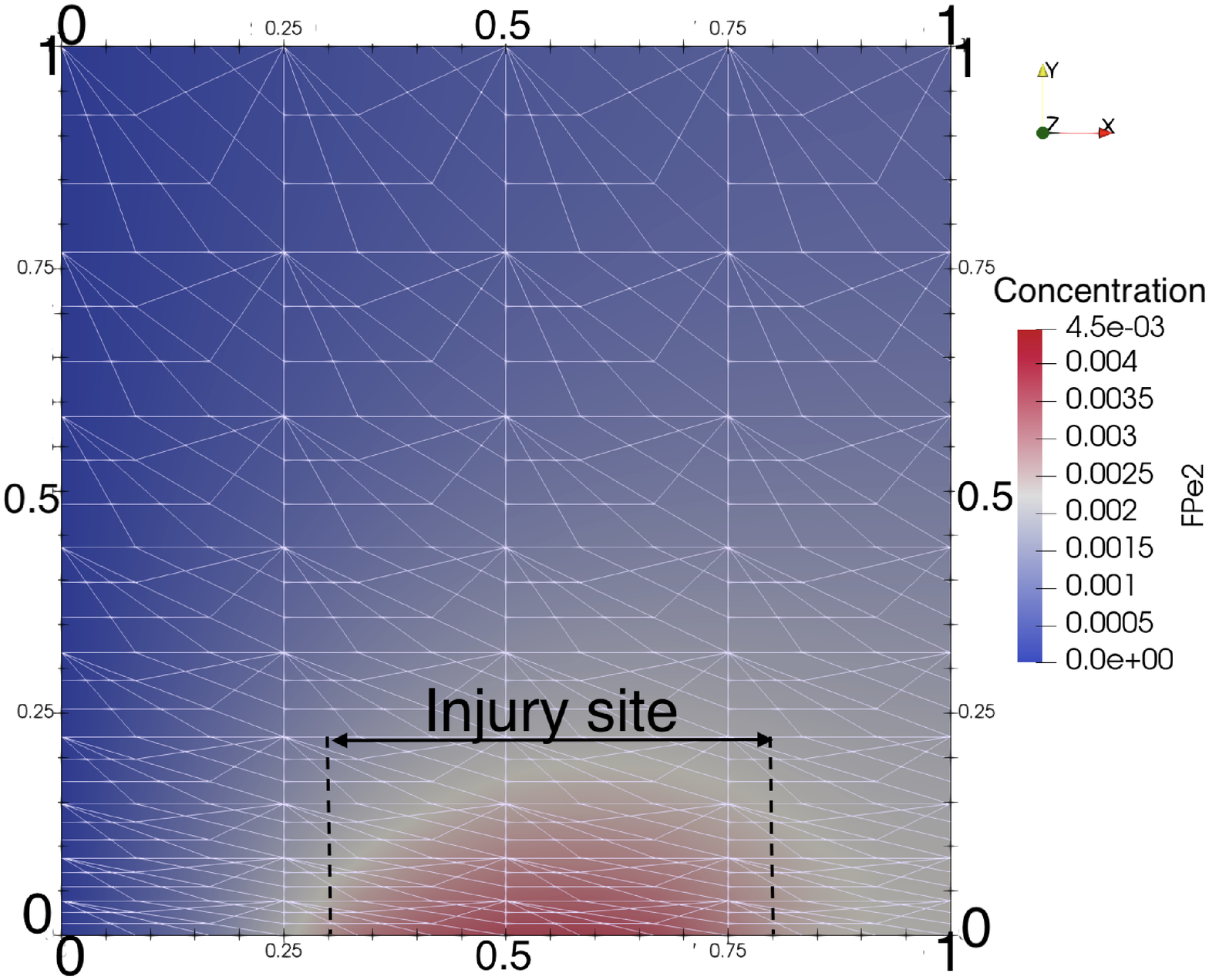}
    \includegraphics[width=0.49\textwidth]{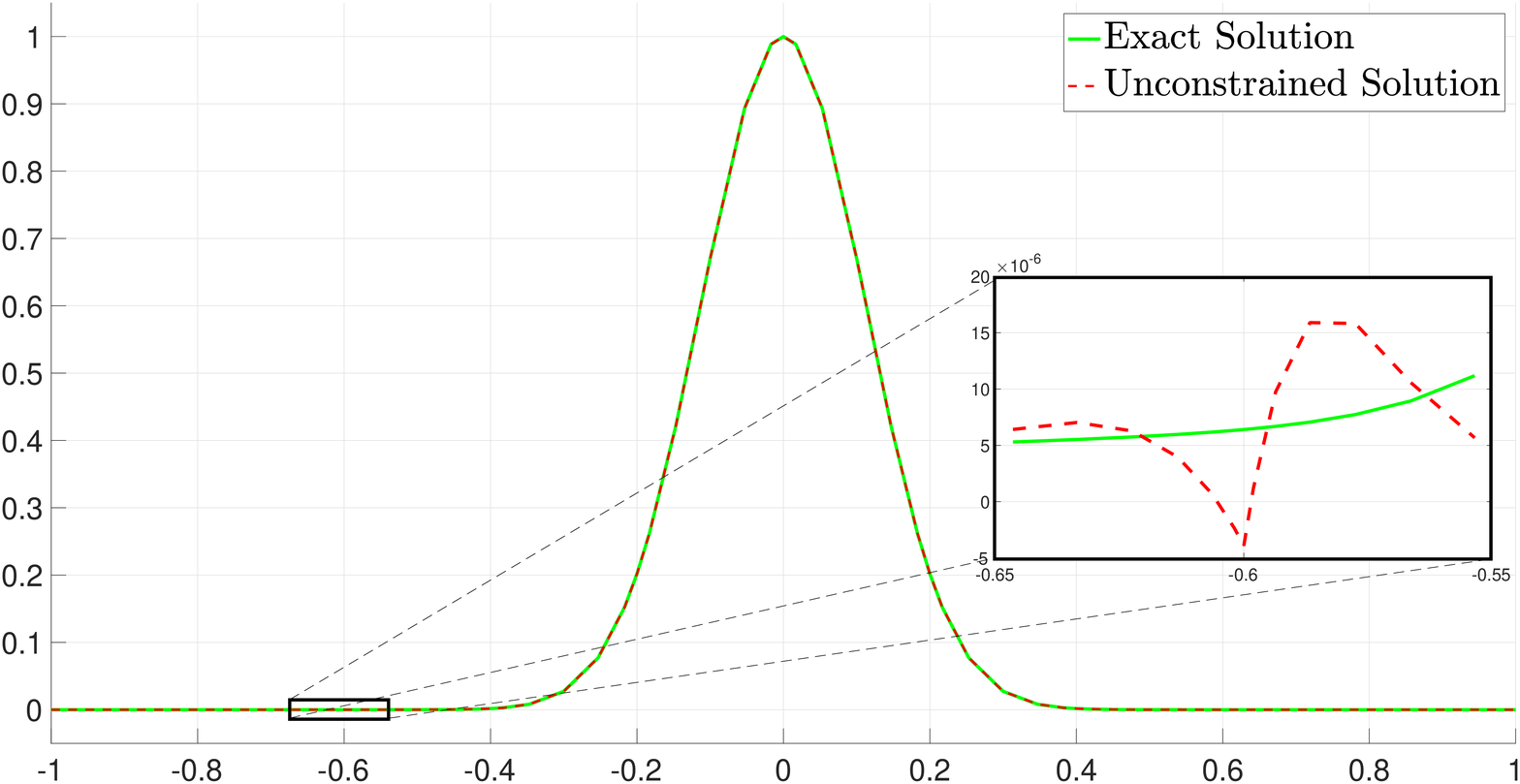}
 \caption{Left: Snapshot in time of a finite element method solution for the platelet aggregation and blood coagulation model \cite{leiderman2011grow} that shows the evolution of the fluid phase chemical Thrombin (FPe2) near the injury site of the vessel. The value of the concentration of FPe2 at a particular point in time is observed to be negative because of the concentration profile near the injury site, which, in turn, causes the simulation to fail.  Right: A simple implementation of the sharp change in concentration in a combustion-like scenario using polynomial projection on a one-dimensional domain. Inset shows the concentration at point $x=-0.6$ is negative (not physically meaningful), which can lead to failure of the simulation.}\label{fig:introeg1}
\end{figure}

\Cref{fig:introeg1} shows examples depicting modeling issues resulting from a nontrivial application. We see that the phenomenon of invalid solutions occurs frequently in many fields, which employ polynomial-based methods for simulations. These simulations are prone to instability or failure as a result of the feedback of an invalid solution from one step to another. It is therefore essential to carefully address the issue of invalid solutions in a generic, domain independent and robust manner without a loss of stability or convergence of the numerical solutions.

\subsection{Contribution}

In this paper, we present a general framework for preserving structure in piecewise polynomial-based time-dependent PDE solvers.  {The procedure, which is agnostic to timestepper by design, is applied at the end of each time-step.} %projections used to numerically solve systems derived out of real-world models. 
We propose a nonintrusive procedure to address the structure-preservation problem. We apply a postprocessing optimization procedure after time steps in the numerical PDE solver that enforces structure as a convex constraint. For certain constraints this optimization is norm-contractive so that our procedure can be interpreted as an application of a (nonlinear) \textit{filter}.

Although the solution we propose corresponds to a conceptually simple convex feasibility problem, the convex feasible set is ``complicated”; in particular, the feasible set is not a finite intersection of simple convex sets, such as Euclidean halfspaces. Therefore, many standard convex optimization algorithms cannot be used to directly solve this problem. However, recent work in \cite{paper0} proposes algorithms to solve this optimization problem for the approximation problem. Our contribution in this paper is to investigate and extend this procedure applied to the numerical solution of PDEs. Our investigations require us to make additional algorithmic advances that are of interest to some PDE solutions: conservation of mass and preservation of boundary conditions and interelement fluxes.  We investigate the efficacy of our filter to efficiently and accurately compute solutions to PDEs while preserving physical structure. All our investigations in this paper are limited to time-dependent PDEs in one spatial dimension.{ We explore a combination of implicit and explicit approaches to solving the PDEs.

The outline of this paper is as follows. In \Cref{sec:survey}, we survey existing solutions to the problem considered in this paper and briefly discuss their limitations. A brief, high-level description of the proposed optimization procedure/filter is given in \Cref{sec:summ}, but a more detailed description of the algorithm including its application to function approximation can be found in \cite{paper0}. In \Cref{sec:summ}, we also summarize notation used throughout the paper and describe the types of constraints we consider. 
%The extension of this concept to PDEs in 1D is the subject matter of this work.  
\Cref{sec:method} details the application to PDEs and some additional discussion of the filtering procedure using geometrical interpretations. We also discuss how the filtering procedure changes some quantities of interest, such as values on element boundaries and total mass. We subsequently formulate a remedy to conserve these quantities by incorporating extra constraints into the filter. \Cref{sec:algo} summarizes the proposed algorithm as formulated in  \Cref{sec:method}. Finally, \Cref{sec:results} presents numerical examples with solutions to PDEs using both discontinuous Galerkin (DG) and continuous Galerkin (CG) formulations. We provide empirical evidence that this optimization procedure has negligible impact on convergence rates and absolute accuracy for these solvers.

\subsection{Existing techniques for structure preservation}\label{sec:survey}

A number of strategies have been proposed in the literature to preserve structure of polynomial approximations.
%while conserving the quantities of interest and avoid invalid solutions. 
Many of these strategies can successfully guarantee the preservation of structure in special cases, or with special discretizations. For narrative purposes, we partition the existing methods into two broad categories: nonintrusive and intrusive.

{Methods are considered to be nonintrusive if they constrain the solution obtained from the solver with minimal, often superficial, change to the numerical scheme. 
Prominent nonintrusive methods include limiters applied at each time-step \cite{TVD,A,B,C,E}. These limiters typically affect only the design of numerical fluxes and not the underlying scheme. However, they must be specially designed for different kinds of constraints and approximation spaces.  These limiters therefore lack flexibility with respect to discretization, spatial dimension, and type of structure/constraint. %This limits the application to the specific problem, domain, dimensionality and function spaces and therefor lack robustness. . 

We also note that many different types of structure-preserving constraints can be considered. For example, the technique from \cite{A,A4} can successfully impose maximum principle-based constraints on a scalar or vector field.  The MVMT-OBR approach described in \cite{a1} prescribes mass and maximum principle conservation successfully on a discrete set of points in the domain.  However, some applications require that although two fields $u$ and $v$ need not individually obey a maximum principle, the sum $u + v$ must obey such a principle \cite{A4}. In such scenarios, standard maximum principle approaches cannot be employed.

The solution we propose belongs to this nonintrusive category, but attempts to mitigate the previously mentioned flexibility issues. We consider a continuous version of the constraint satisfaction problem that ensures the constraints are satisfied on all the points in the domain, and not on just a discrete set of points. Furthermore, the proposed method can be used for any arbitrarily high polynomial orders without changes to the algorithm. 

{The second class of methods are those that are intrusive. An approach is considered to be intrusive if it needs to substantially change the underlying numerical scheme or properties of the solution domain.} The intrusive methods look at the structure-preservation problem as a PDE-constrained optimization problem. Some of these methods modify the spatial discretization \cite{zahr2018optimization}, and others use limiters derived from Karush-Kahan-Tucker (KKT) optimality conditions \cite{KKT}.  The strategies proposed for constraint satisfaction in \cite{a1} consider a version of the problem that is a subset of the one solved by the proposed solution.  In \cite{a2} and its extensions \cite{a1,l_1}, the authors explore approaches for positivity preservation that are based on basis functions derived from Bernstein polynomials such that they are non-negative and possess partition of unity property.  Therefore the interpolated solution respects the original bounds at any point in the domain.
%using KKT limiters wherein the constraints are folded as part of the PDE system and the subsequent numerical strategy uses this modified system of PDE. 
While successful in imposing structure as part of the PDE discretization, intrusive approaches often suffer from the limitation of having to modify the numerical scheme, and the type of modification is both problem- and constraint-dependent. It requires substantial human intervention to change the discrete PDE solver. Methods involving changes in domain to solve this problem are also largely problem-specific and therefore lack flexibility. For example, the authors in \cite{zahr2018optimization} use an optimization problem incorporated in the scheme, and the solution is computed on a curved mesh that tracks discrepancies. 
 
To summarize, existing strategies in the literature to preserve structure typically come in the form of intrusive methods, requiring nontrivial modification to numerical schemes, or nonintrusive methods, which typically affect existing numerical implementations in benign ways.The procedure we consider in this paper falls into the latter category, and our framework handles very general constraints. Furthermore, the mathematical formulation is agnostic to the type of (linear, convex) constraint and the spatial dimension of the problem. The price we pay for this generality is that some nontrivial (but convex) optimization must be performed. We describe this optimization in more detail in the next section.

\section{Structure-preserving function approximation} \label{sec:summ}

In this section, we summarize the main algorithmic ideas from \cite{paper0}, which is a major ingredient for our approach. This approach is a map $M$ from a given function $u$ (that may or may not satisfy structural constraints) in a finite-dimensional space $V$  to a unique function $M(u) \in V$, which does satisfy these constraints. The work in \cite{paper0} defines a general class of constraints, corresponding to a feasible set in $V$ that is an affine convex cone, and shows that $M(u)$ is the projection of $u$  onto this feasible set, which is an affine convex cone in $V$. Unfortunately, this cone is not a polytope, and convenient parameterizations of $V$ result in representation as an intersection of an (uncountably) infinite number of halfspaces. This formulation and parameterization does not easily lend itself to existing algorithms, so \cite{paper0} develops some novel algorithms, based on seminal convex feasibility algorithms \cite{von_neumann_functional_1951}. We now briefly discuss the types of constraints considered in \cite{paper0} and some algorithms to implement them.

Let $\Omega \subset \R^d$ be a physical domain. In this paper, we are interested in $d = 1$, but this restriction is not necessary for the general approach. Let $V$ be a finite-dimensional Hilbert space of real-valued functions on $\Omega$ (for example, polynomials up to some fixed, finite degree), and suppose $u \in V$ is a given function. The approach in \cite{paper0} considers families of constraints, each of the form
\begin{align}\label{eq:constraints}
  \mathcal{L}_x(u) &\leq \ell(x), & x &\in \Omega,
\end{align}
where $\mathcal{L}_x$ is a linear operator that is bounded on $V$, and $\ell$ is a function on $\Omega$. The feasible set in $V$ adhering to such a constraint family is convex and includes the following examples:
\begin{itemize}
\item Positivity: $u(x) \ge 0$ for all $x$ in $\Omega$.
\item Monotonicity: $u'(x) \ge 0$ for all $x$ in $\Omega$.
%\item Boundedness: $0\le p(x) \le 1 $ for all $x$ in $ \Omega_e$.
\end{itemize}
Note that the framework in \cite{paper0} allows a finite number of such families to be considered simultaneously, so that boundedness, e.g., enforcing $0 \leq u(x) \leq 1$, is also a valid constraint. The families of constraints correspond to a feasible set {$K \subset V$} corresponding to the elements of $V$ that satisfy the constraints. The strategy in \cite{paper0} is to solve the convex feasibility problem
\begin{align}\label{eq:V-opt}
  M(u) \coloneqq \argmin_{k \in K} \| u - k\|_V,
\end{align}
which is well posed. If {$0 \in K$}, then this optimization problem is norm-contractive \cite[Proposition 5.1]{paper0}, and therefore can be interpreted as a filter.

We define more notations to describe the algorithm. Suppose $V$ is an $N-$dimensional subspace of a Hilbert space $H$, with $\{\psi_j\}_{j=0}^{N-1}$ a collection of orthonormal basis functions,
\begin{align*}
  V &= \mathrm{span}\left\{ \psi_0, \ldots, \psi_{N-1} \right\}, & \left\langle \psi_i, \psi_j \right\rangle &= \delta_{ij}, & i, j &= 0, \ldots, N-1,
\end{align*}
where $\left\langle\cdot,\cdot\right\rangle$ is the inner product on $V$, and {$\delta_{ij}$}, the Kronecker delta function.  { For a particular constraint {$k \in K$}, we can represent $u \in V$ in its coordinates $\{\hat{v}_j\}_{j=0}^{N-1}$ collected in a vector $\bs{{v}} \in \R^N$}. {Any $u \in V$ that does not satisfy the desired constraints can be represented as $\sum\limits_{j = 0}^{N-1} \tilde{v}_j \psi_j =   \bs{\tilde{v}}\bs{\psi}$. We collect the coefficients of expansion in a vector $\bs{\tilde{v}} \in \R^N$.} The optimization problem \eqref{eq:V-opt} is therefore equivalent to
{
\begin{align}\label{eq:v-optimization}
  \argmin_{\bs{{v}} \in C} \| \tilde{\bs{v}} - \bs{{v}}\|_2,
\end{align}}
{where $\bs{v}$ is the filtered version of $ \tilde{\bs{v}}$ and obeys the constraints, }$\|\cdot\|_2$ is the Euclidean 2-norm on vectors, and $C$ is the affine conic region in $\R^N$ corresponding to {$K \subset V$}.  Whereas the basis function $\bs{\psi}$ represents any orthonormal basis function, in general, use of any basis function is generally possible as long as a transformation to the orthonormal basis is done prior to the application of the filter. Note that the filter operates on the coefficients of expansion ${\tilde{\bs{v}}}$ while ensuring that the filtered $\bs{v}$  can be mapped back to the baseline set of basis function, thus maintaining the constraint satisfaction on all the points in the domain. If {$K$} contains a single family of the form \eqref{eq:constraints}, then the set $C$ can be written as
\begin{align*}
  C = \bigcap_{x \in \Omega} H_x,
\end{align*}
where $H_x$ are halfspaces in $\R^N$. In other words, for a fixed $x$, $H_x$ is an $(N-1)$-dimensional planar surface in $\R^N$ defined by the single linear constraint \eqref{eq:constraints}. The algorithms in \cite{paper0} proceed by computationally inspecting the signed distance function,
\begin{align*}\label{eq:distform}
  s(x) \coloneqq \mathrm{sdist}(\bs{\tilde{v}}, C_x) = \left\{ \begin{array}{cc} -\mathrm{dist}(\bs{\tilde{v}}, H_x), & x \not \in C_x, \\
                                                          +\mathrm{dist}(\bs{\tilde{v}}, H_x), & x \in C_x.
  \end{array}\right.
\end{align*}
Here, ``inspection" means, for example, the ability to compute the global minimum of $s(x)$ and/or to determine regions where $s$ is negative. Based on this inspection, the algorithms project the state vector of the current iterate $\bs{\tilde{v}}$ onto $H_y$ for some $y \in \Omega$, or perform relaxed/averaged versions of these projections. Geometrically, the projection operation corresponds to projecting $\bs{\tilde{v}}$ onto a supporting hyperplane $H_x$ for $C$, which we depict in \Cref{fig:distandprojection}.

\begin{figure}[htbp]
  \begin{center}
    \resizebox{\textwidth}{!}{
        \includegraphics[width=0.49\textwidth]{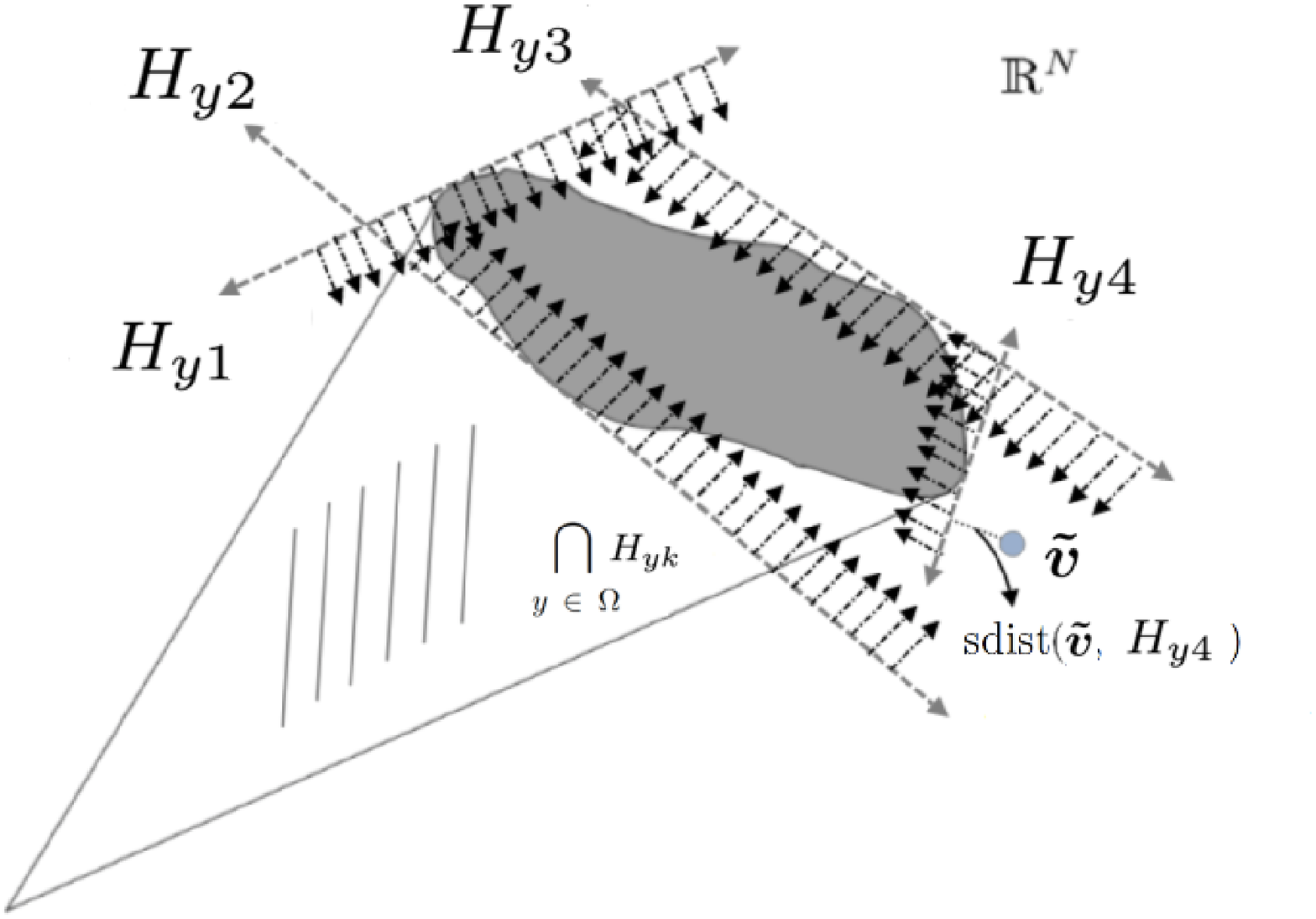}
      \includegraphics[width=0.49\textwidth]{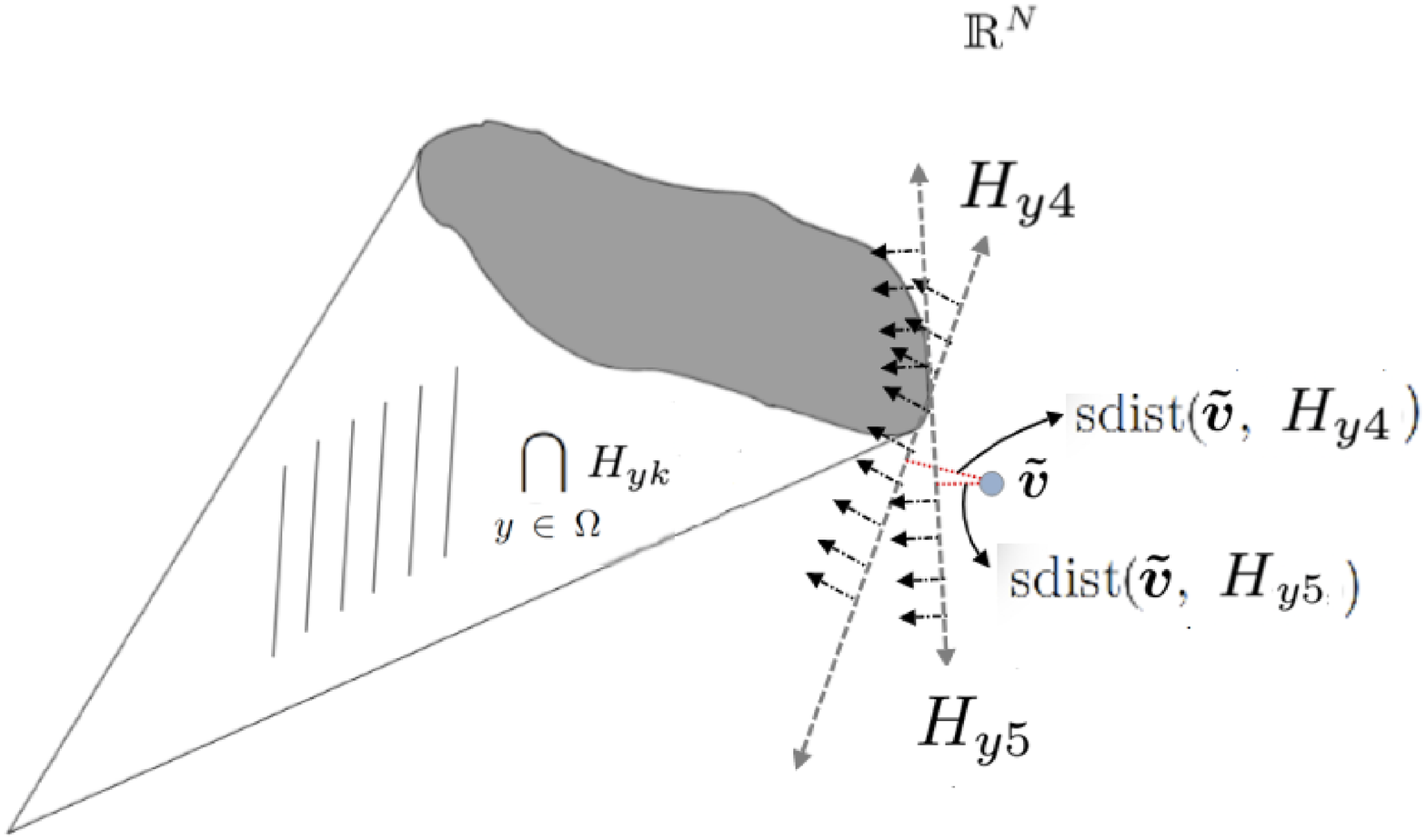}
    }
  \end{center}
 \caption{The steps in the procedure that find the distance between $\bs{\tilde{v}}$ {and the hyperspace boundary representing} each violating constraint. Subsequently, the algorithm greedily calculates the correction to $\bs{\tilde{v}}$.
 Left: A geometrical visual of the distance calculation from $\bs{\tilde{v}}$ to the hyperplanes defining boundaries of the violating constraint.  Right: Selection of $H_{y4}$ over $H_{y5}$ since it defines the hyperplane farther away compared to other violating constraints. Projection of $\bs{\tilde{v}}$ on to $H_{y4}$.}\label{fig:distandprojection}  
 \end{figure}

The process is repeated until numerical convergence up to a tolerance (i.e., until the minimum value of the signed distance function is numerically 0). Thus, this algorithm is a type of generalized iterative cyclic/alternating algorithm applied to the case of an infinite number of convex sets (halfspaces). 

The major computational work in this iterative procedure is the manipulation/minimization of $s(x)$ at each iteration. This signed distance function has the form
\begin{align}\label{eq:sdist}
  s(x) &= \lambda(x) \left( \mathcal{L}_x(u) - \ell(x) \right), & \lambda^2(x) \coloneqq \frac{1}{\sum_{j=0}^{N-1} \left(\mathcal{L}_x(\psi_j)\right)^2},
\end{align}
where, $\ell \in \Omega$, and $s$ is a $\lambda$-weighted version of $u$, therefore it is easy to evaluate. However computing a global minimum for $s$ may be difficult in general. In \cite{paper0}, it is shown that if $V$ is a univariate polynomial space of degree $N-1$, then computing the minimum of $s$ can be accomplished by computing the roots of a polynomial of degree $3 N$. We accomplish this by computing the spectrum of a confederate matrix associated with a Legendre orthogonal polynomial basis. For more details on the algorithm, including computational cost, see \cite{paper0}.

This algorithm which seeks to solve the optimization problem \eqref{eq:V-opt} is the key ingredient in our approach to preserve structure in time-dependent PDE simulations.

%% file: Method.tex
\section{Method}\label{sec:method}

We now discuss the application of the filter introduced in \Cref{sec:summ} to time-dependent PDEs in one spatial dimension. We will primarily focus on method-of-lines discretizations with a Galerkin-type spatial discretization. In particular, we will consider continuous and discontinuous Galerkin formulations. Consider an advection-diffusion-reaction system defined as

\begin{align}\label{eq:contadr}
{u}_t(x,t) + a \cdot {u}_x(x,t) = \gamma {u}_{xx}(x,t) + {r}({u}(x,t)) 
\end{align} 

\noindent where $x \in \Omega$, $a$ is the velocity for advection, $\gamma$ is the diffusivity, and ${r}$ is a non-linear function representing the reaction term. For simplicity, assume $a$ and $\gamma$ to be constant advection and diffusion coefficients, respectively. Assuming appropriate initial $\bs{u}_0(x,t=0)$ and boundary conditions are defined for \eqref{eq:contadr}, we can formulate its semidiscrete form as follows. 

Galerkin-type methods assume an ansatz for $u$ as a time-varying element of a fixed $N$-dimensional linear subspace $V$, where frequently $V \subset L^2(\Omega)$:
\begin{align}\label{eq:u-ansatz}
  u(x,t) \approx u_N(x,t) &\coloneqq \sum_{i=0}^{N-1} \tilde{v}_i(t) \phi_i(x), & V &= \mathrm{span} \{ \phi_1, \ldots, \phi_N \}.
\end{align}
Note that the basis functions $\phi_i$ used in \eqref{eq:u-ansatz} represent the traditional FEM basis functions that span the entire $\Omega$ and are not necessarily orthogonal. We denote the orthonormal basis functions by $\bs \psi_j$. We assume that both basis sets span the same space and thus a transformation between them exists.  For example, finite element methods partition $\Omega$ into {$E$} nonoverlapping subintervals.  As a next step, the method assumes that the continuous functions in $V$ are polynomials of a fixed degree $N$ on each $E \in \Omega$. The discontinuities in derivatives are allowed only at partition boundaries. Similarly, discontinuous Galerkin finite element methods define $V$ in a similar way except that elements are allowed to be discontinuous at partition boundaries. The semidiscrete form for \eqref{eq:contadr} is derived in the standard Galerkin way, by using the ansatz \eqref{eq:u-ansatz} and forcing the residual to be $L^2$-orthogonal to $V$. Usually, integration by parts is performed in the residual orthogonalization step and often, depending on the equation and spatial discretization, a  numerical flux and/or stabilization terms are included in the resulting weak formulation.

The result is a system of ordinary differential equations prescribing time-evolution of the discrete degrees of freedom represented by vector $\bs{\tilde{v}} =  \{\tilde{v}_0,\ldots,\tilde{v}_{N-1}\}$.

\begin{align}\label{eq:disadr}
  \bs{M} \frac{\partial}{\partial t}{\bs{\tilde{v}}} + \bs{A}{\bs{\tilde{v}}} = -\gamma \bs{L} {\bs{\tilde{v}}}  + {\bs{r}} + \bs{F}(\bs{\tilde{v}})
\end{align}
where $\bs{M}, \bs{A}$, and $\bs{L}$ are the $N \times N$ mass, advection, and Laplacian (stiffness) matrices, respectively, defined as
\begin{align*}
  (M)_{i,j} &= \left\langle \phi_i, \phi_j \right\rangle, & 
  (A)_{i,j} &= \left\langle \phi_i, a \frac{\partial}{\partial x} \phi_j(x) \right\rangle, & 
  (L)_{i,j} &= \left\langle \frac{\partial}{\partial x} \phi_i, \frac{\partial}{\partial x} \phi_j \right\rangle.
\end{align*}
The $N$-vector ${\bs{r}}$, defined as
\begin{align*}
\bs{r}= \sum_{j=0}^{N-1} \hat{r}_i(t) \phi_i(x), & V &= \mathrm{span} \{ \phi_1, \ldots, \phi_N \},
\end{align*}
has entries that are numerical approximations to the integral of the nonlinear reaction term,
\begin{equation*}
  \hat{r}_i(\bs{\tilde{v}},t) \approx \int \bs{r}(u_N,t) \phi_i(x) \dx{x},
\end{equation*}
which in this paper we compute with a collocation-based approach. Finally, the term $\bs{F}(\bs{\tilde{v}})$ is a {generic} term for any numerical fluxes or stabilization terms. For example, in an advection-dominated problem with a DG formulation, $\bs{F}$ might be the upwind flux corresponding to the continuous advection term $a \bs{u}_x$. Finally, a fully discrete scheme is derived from \eqref{eq:disadr} using an appropriate time-integration method. 

Let $\bs{{v}}^n$ represent the solution of \eqref{eq:disadr} at timestep $n$, and let us call $\bs{\tilde{v}}^{n+1}$ the solution at the timestep $n+1$.
For simplicity, assume that the solutions $\bs{{v}}^n$ and $\bs{\tilde{v}}^{n+1}$ are transformed into orthonormal basis $\bs{\psi_j}$ and back to original basis $\phi_I$ inside the filter as the first and last steps. This discrete scheme does not enforce the structural properties that we desire. As discussed in the \Cref{sec:summ}, given a constraint set $C \subset \R^N$, if $\bs{v}^n \in C$, then it is not necessarily true that $\bs{\tilde{v}}^{n+1} \in C$. To rectify this situation, we employ the optimization outlined in  \Cref{sec:summ} as a postprocessing step. We will hereafter refer to this optimization as a nonlinear ``filter" due its norm-contractivity properties for the types of constraints we consider \cite[Proposition 5.1]{paper0}. Thus, our proposed procedure is a simple, {nonintrusive} augmentation of the standard fully discrete scheme \eqref{eq:disadr}:
\begin{equation}\label{filter}
  \bs{{v}}^n \xrightarrow{\textrm{Timestepper for \eqref{eq:disadr}}}  \bs{\tilde{v}}^{n+1} \xrightarrow{\textrm{Filter}} \bs{{v}}^{n+1}
\end{equation}  

\begin{figure}[H]
  \centering
\includegraphics[width=0.6\textwidth]{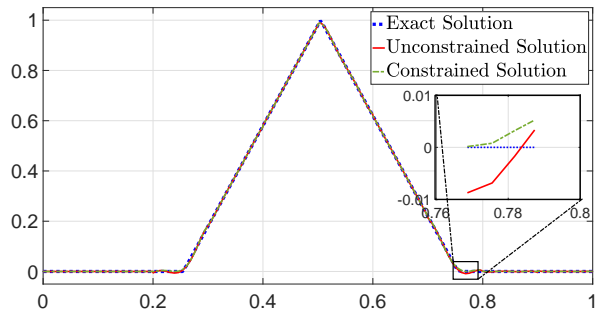}   
    \caption{Numerical solution to a PDE with and without application of the positivity-constraining filter. Shown are the exact solution, the unconstrained solution (i.e., the solution via the scheme \eqref{eq:disadr}), and the constrained solution (the solution via \eqref{filter}).} \label{fig:1} 
\end{figure}

Thus, we obtain the filtered solution $\bs{{v}}^{n+1} = \sum_{j = 0}^{N-1} \hat{v}_j(t) \phi_j(x)$ and introduce it  as the state in the next step of the fully discrete scheme. The effect of the filter for an example of constraint: positivity preservation for $\bs{\tilde{{v}}}^n $ is demonstrated in \Cref{fig:1}, which shows a snapshot in time of an advecting wave. The property of interest here is the positivity, and in this particular example, a discontinuous Galerkin formulation is used as the PDE solver. The process of enforcing positivity (see inset) leads to changes in solution properties; in particular, elementwise boundary values and the mean of the discrete solution are not preserved. A change in element boundary values also changes the corresponding (say, upwind) fluxes at the next timestep. Hence, in this particular case, the filter changes physical behavior (numerical fluxes and mean value). 
 
We will assume hereafter that numerical fluxes are computed as explicit functions of left- and right-boundary values. Hence enforcing element boundary value conservation for a discontinuous Galerkin solver becomes necessary to ensure that numerical fluxes at future times are unchanged by this procedure. One remedy for this particular issue is to impose additional equality constraints (i.e., function values at element boundaries) in the filter. The preservation of the mean value (e.g., total mass) can also be enforced in a similar way. Later in this section we elaborate on the construction of the filter for these additional constraints. Note that, between CG and DG discretizations, element boundary value constraints are relevant (for the purposes of conserving numerical fluxes) only for DG-type discretizations.

\subsection{Degrees of freedom in CG and DG simulations}\label{ssec:dof}

In order to describe our enforcement of equality constraints, we need to introduce some new notation that is specific to the physical discretization. Since DG discretizations have degrees of freedom that are decoupled across elements, the $N$ degrees of freedom in the optimization described in \Cref{sec:summ} only need to be over element-local degrees of freedom. In constrast, standard degrees of freedom in CG simulations (i.e., coefficients of ``hat" and ``bubble" functions) are coupled in $L^2$ across elements. 

\subsubsection{Discontinuous Galerkin}\label{sssec:dg}
With a partition consisting of $E$ subintervals of $\Omega$, a DG discretization allows us to convert the optimization \eqref{eq:v-optimization} in an $N$-dimensional space into a much more efficient and parallelizable set of $E$ independent optimizations each in $n \coloneqq N/E \ll N$ dimensions. First, we recognize that a DG expansion of the form \eqref{eq:u-ansatz} can be written as 
\begin{align*}
  u_N(x,t) &= \sum_{e=1}^E \sum_{k=0}^{n-1} \tilde{v}_{e,k} \psi_{e,k}(x), & N &= E \, n
\end{align*}
where $\{\psi_{e,k}(\cdot)\}_{k=0}^{n-1}$ for a fixed $e$ are polynomials whose support is \emph{only} on element $e$. In particular, these element-local polynomials can be chosen as $L^2(\Omega)$-orthonormal (e.g., mapped Legendre) polynomials. Therefore, 
\begin{subequations}\label{eq:dg-assumptions}
\begin{align}
  \left\langle \psi_{e,k}, \psi_{f,\ell} \right\rangle &= 0, & e \neq f, &\;\;\;\;\; k, \ell \in \{0, \ldots n-1\}.
\end{align}
The second observation we make is that the linear constraint operator $\mathcal{L}_x$ defined in \eqref{eq:constraints}, for common examples we consider, also obeys this type of decoupling. If $\{\Omega_e\}_{e=1}^E$ is the partition of $\Omega$, then we assume that
\begin{align}
  \mathcal{L}_{x}(\psi_{e,k}(x)) &= 0, & x &\not\in \Omega_e.
\end{align}
\end{subequations}
This assumption is satisfied for point and derivative evaluation, i.e., for $\mathcal{L}_x(u) = u(x)$ and $\mathcal{L}_x(u) = u'(x)$. Under conditions \eqref{eq:dg-assumptions}, both the constraint and the $L^2$ norm are decoupled across elements, and so we have the following result:
\begin{proposition}\label{prop:dg-decouple}
  Assume that conditions \eqref{eq:dg-assumptions} hold. Given any $\tilde{\bs{v}} \in \R^N$, partition it into $E$ $n$-dimensional subvectors, each containing element-local degrees of freedom:
  \begin{align}\label{eq:dg-v-partition}
    \tilde{\bs{v}} &= \left(\begin{array}{c} \tilde{\bs{v}}_1 \\ \tilde{\bs{v}}_2 \\ \vdots \\ \tilde{\bs{v}}_E \end{array}\right), & 
      \tilde{\bs{v}}_e &= \left(\begin{array}{c} \tilde{v}_{e,0} \\ \vdots \\ \tilde{v}_{e,n-1} \end{array}\right).
  \end{align}
  Then the solution to \eqref{eq:v-optimization} is given by
  \begin{align*}
    \argmin_{\bs{v} \in C} \left\| \bs{v} - \tilde{\bs{v}} \right\|_2 = \left(\begin{array}{c} 
      \argmin_{\bs{w} \in C^1} \left\| \bs{w} - \tilde{\bs{v}}_1 \right\|_2 \\
      \argmin_{\bs{w} \in C^2} \left\| \bs{w} - \tilde{\bs{v}}_2 \right\|_2 \\
      \vdots \\
      \argmin_{\bs{w} \in C^E} \left\| \bs{w} - \tilde{\bs{v}}_E \right\|_2
    \end{array}\right),
  \end{align*}
  where the $n$-dimensional sets $C^e$ are defined as
  \begin{align*}
    C^e \coloneqq \left\{ \bs{w} = (w_0, \ldots, w_{n-1})^T \in \R^n \; \big|\; \sum_{k=0}^{n-1} w_k \psi_{e,k}(x) \leq r(x) \; \forall\; x \in \Omega_e \right\}.
  \end{align*}
\end{proposition}
We emphasize again that DG discretizations satisfy \eqref{eq:dg-assumptions} so that Propostion \ref{prop:dg-decouple} allows us to conclude that our optimization (for example constraining positivity and/or monotonicity) can be decoupled to operate on individual elements. The decoupling can result in substantial computational savings since frequently one $N$-dimensional optimization is much more expensive than $E$ $N/E$-dimensional optimizations.
\begin{remark}
  The elementwise decoupling conclusion of Proposition \ref{prop:dg-decouple} holds under slightly more general conditions. For example, in so-called $p$-adaptive simulations, the number of local degrees of freedom may depend on the element index, i.e., $n = n(e)$. However, we do not pursue $p$-adaptive simulations in this paper, so this level of generality is not needed.
\end{remark}

\subsubsection{Continuous Galerkin}
In contrast to the discontinuous Galerkin framework, in continuous Galerkin discretizations, typically both the constraint operator (e.g., point evaluation operator) and the $L^2$ norm are coupled across elements. Thus, there is no straightforward decoupling procedure that can be leveraged. Instead, the full optimization problem \eqref{eq:v-optimization} in $N$-dimensional space must be solved. 

An additional difficulty with this optimization in CG simulations is that the discussion and algorithms presented in Section \ref{sec:summ} essentially require an expansion in \emph{orthonormal} basis functions, i.e., \eqref{eq:u-ansatz} must be expressed as
\begin{align}\label{eq:cg-orthonormal}
  u_N(x,t) &= \sum_{j=0}^{N-1} \tilde{v}_j(t) \phi_j(x) = \sum_{j=0}^{N-1} \tilde{w}_{j}(t) \psi_j(x), &
  \left\langle \psi_j, \psi_k \right\rangle &= \delta_{k,j}.
\end{align}
Typically, CG simulations utilize non-orthonormal hat and bubble functions as degrees of freedom \cite{karniadakis2013spectral} an therefore the use of this optimizer requires transformation of the hat-bubble coordinates into an orthonormal set of coordinates. To accomplish this, (a Cholesky factor of) the mass matrix $\bs{M}$ must be inverted. Typically $\bs{M}$ is sparse, but its inverse Cholesky factor is usually not, which poses a challenge for large-$N$ simulations when size-$N$ dense matrix linear algebra is computationally infeasible.

In this paper, we consider simulations only in one spatial dimension where $N$ is small enough so that direct inversion of $\bs{M}$ is feasible. However, more sophisticated procedures would be needed to extend this approach to two or three spatial dimensions.

\subsection{Elementwise boundary value conservation (DG)}\label{sec:fluxcons}

The goal in this section is to enforce the element boundary constraints explained at the end of \Cref{sec:method} for DG discretizations. For simplicity of exposition, we assume in what follows that $\{\psi_{e,k}\}_{k=0}^{n-1}$ are $L^2$-orthonormal. The procedure we describe can be extended to the case when this assumption is not satisfied. With $\{\Omega_e\}_{e=1}^E$ the subinterval partition of $\Omega$, we seek to solve the equality-constrained optimization problem,
\begin{align}\label{eq:eq-const-full}
  \min_{\bs{v} \in C} \left\| \bs{v} - \tilde{\bs{v}} \right\|_2 \hskip 5pt \textrm{such that} \hskip 5pt \sum_{k=0}^{n-1} \tilde{v}_{e,k} \psi_{e,k}(x_e^\pm) = \sum_{k=0}^{n-1} v_{e,k} \psi_{e,k}(x_e^\pm) \hskip 5pt \forall \; e \in \{1, \ldots, E\},
\end{align}
where $\tilde{v}_{e,k}$ and $v_{e,k}$ are components of the elementwise subvector partition of $\tilde{\bs{v}}$ and $\bs{v}$, respectively, that was introduced in \eqref{eq:dg-v-partition}. The point values $x_e^{\pm}$ are left- and right-hand side values of subinterval $e$, $\Omega_e = [ x_e^-, x_e^+ ]$. As described in \Cref{sssec:dg} and Proposition \ref{prop:dg-decouple}, the (non-equality-constrained) DG optimization problem can be decoupled into elementwise operations. Adding in the elementwise constraints described above does not change this decoupling property (since the boundary values on element $e$ are independent of the degrees of freedom on any other element). The resulting optimization on element $e$ has the form,
\begin{align}\label{eq:eqconst-opt}
  \min_{\bs{w} \in C^e} \left\| \bs{w} - \tilde{\bs{v}_e} \right\|_2 \hskip 5pt \textrm{such that} \hskip 5pt \bs{Q}^T \left( \bs{w} - \tilde{\bs{v}}_e \right) = \bs{0},
%  \tilde{v}_{e,k} \psi_{e,k}(x_e^\pm) = \sum_{k=0}^{n-1} v_{e,k} \psi_{e,k}(x_e^\pm).
\end{align}
where $\bs{Q}$ is an $n \times 2$ matrix with orthonormal columns $\bs{q}_1, \bs{q}_2$ satisfying
\begin{align*}
  \mathrm{span}\left\{ \bs{q}_1, \bs{q}_2 \right\} &= \mathrm{span}\left\{ \bs{\psi}_e(x_e^-), \bs{\psi}_e(x_e^+) \right\}, & \bs{\psi}_{e}(x) &\coloneqq \left( \psi_{e,0}(x), \; \cdots \; \psi_{e,n-1}(x)\right)^T.
\end{align*}
In order to solve this $n$-dimensional linear-equality-constrained optimization problem, we reduce it to an $(n-2)$-dimensional optimization problem of the standard form \eqref{eq:v-optimization} by working in the $n-2$ coordinates corresponding to $\mathcal{R}(\bs{Q})^\perp$. To set up for this procedure, let $\{\bs{p}_1, \bs{p}_2, \ldots, \bs{p}_{n-2}\}$ be an(y) orthonormal completion of $\{\bs{q}_1, \bs{q}_2\}$ in $\R^n$, and introduce $\bs{P} \in \R^{n \times (n-2)}$,
\begin{align*}
  \bs{P} &= \left(\begin{array}{cccc} \bs{p}_1 & \bs{p}_2 & \cdots & \bs{p}_{n-2} \end{array}\right), %& \left\langle \bs{p}_j, \bs{q}_k \right\rangle &= 0, & \left\langle \bs{p}_j, \bs{p}_r \right\rangle &= 0,
\end{align*}
%This matrix is defined such that $(\bs{Q} \; \bs{P}) \in \R^{n \times n}$ is an orthogonal matrix. 
Then, our equality constrained optimization problem \eqref{eq:eqconst-opt} is equivalent to the non equality-constrained problem
\begin{align}\label{eq:eq-reduced}
  \min_{\bs{z} \in \hat{C}} \left\| \bs{z} - \bs{P}^T \tilde{\bs{v}_e} \right\|_2,
\end{align}
where $\hat{C}$ is the $\R^{n-2}$-dimensional set,
\begin{align*}
  \hat{C} &\coloneqq \left\{ \bs{z} \in \R^{n-2} \; \big| \; \sum_{j=0}^{n-3} z_j \hat{\psi}_j(x) \leq \hat{\ell}(x) \right\}, & \hat{\ell}(x) &\coloneqq \ell(x) - \bs{\psi_e}(x)^T \bs{Q} \bs{Q}^T \tilde{\bs{v}_e},
\end{align*}
and $\{\hat{\psi}_j\}_{j=0}^{n-3}$ are the $\bs{P}$-projected functions,
\begin{align*}
  \bs{\psi}_e^T(x) \bs{P} \eqqcolon \left(\begin{array}{ccc} \hat{\psi}_0(x) & \cdots & \hat{\psi}_{n-3}(x) \end{array}\right).x`
\end{align*}
Therefore, the reduced problem \eqref{eq:eq-reduced} is precisely a special case of our problem \eqref{eq:v-optimization}, and all the same tools described in Section \ref{sec:summ} are applicable. The signed distance function in \eqref{eq:sdist} requires knowledge only of the appropriate $n-2$ projected orthonormal functions. The formula precisely is
\begin{align*}
  s(x) &= \lambda(x) \left( \mathcal{L}_x(z) - \hat{\ell}(x) \right), & \lambda^2(x) &= \frac{1}{\sum_{j=0}^{n-3} (\mathcal{L}_x(\hat{\psi}_j(x))^2},
\end{align*}
and $z = \sum_{j=0}^{n-3} z_j \hat{\psi}_j(x)$. Once this solution, say $\bs{z}^\ast$, is computed, then the solution to \eqref{eq:eq-const-full} is the $n$-dimensional vector $\bs{Q} \bs{Q}^T \tilde{\bs{v}}_e + \bs{P} \bs{z}^\ast$.

To summarize, we wish to solve the $N$-dimensional problem \eqref{eq:eq-const-full}, which simultaneously enforces inequality constraints (such as positivity) and preserves element boundary values (which we take as a proxy for preservation of numerical fluxes). This problem is equivalent to solving \eqref{eq:eqconst-opt} for every element index $e$. For a fixed $e$, this latter problem, in turn, is equivalent to the $(n-2)$-dimensional problem \eqref{eq:eq-reduced}, which can be solved with the techniques outlined in Section \ref{sec:summ}. Note in particular that all the matrices involved in this section can be precomputed and stored for use in online time-dependent simulations. Furthermore, if all physical elements are templated on a standard (reference) element, as is common in finite element code, then only one copy of all these size-$n$ matrices is required for the entire simulation.

\subsection{Mass conservation (CG and DG)}\label{sec:energycons}

A similar approach for mass (integral) conservation can be obtained for both CG and DG discretizations using a procedure that is essentially identical to the one described in the previous section. For example, in a CG discretization with orthonormal coordinates $\tilde{\bs{w}}$, we might wish to perform the optimization \eqref{eq:v-optimization} for $\bs{w}$ subject to the equality constraint,
\begin{align*}
  \int_\Omega \sum_{i=0}^{N-1} w_i \psi_i(x)dx = M \coloneqq \int_\Omega \sum_{i=0}^{N-1} \tilde{w}_i \psi_i(x) \dx{x},
\end{align*}
where we recall that $\{\psi_i\}_{j=0}^{N-1}$ are the orthonormal basis functions in \eqref{eq:cg-orthonormal}.  This is equivalent to the equality constraint
\begin{align*}
  \bs{q}^T \bs{w} &= M, & q_j &\coloneqq \int_\Omega \psi_j(x) \dx{x},
\end{align*}
and therefore the procedure of \Cref{sec:fluxcons} can be leveraged (with only one equality constraint and with an $N$-dimensional state vector).  This process leads to the global conservation of mass. In many cases, having local mass conservation and constraint satisfaction across individual elements is desirable. This type of  problem is possible to formulate by careful reconstruction of constraints without any changes to the algorithm. The constraint on the entire domain can be broken down per element, resulting in ``E" constraints, where E = number of elements $ \{e_0, e_1, \cdots,e_{E-1} \} \in \Omega$.  Therefore, the constraint set looks like the one described in \eqref{eq:eq-const-full}. The only difference in the problem formulation will be the definition of the constraint set, which will be a matrix of size $N \times E$, with each column representing a mass conservation condition for each element, respectively. Since we use a linear programming solver in the algorithm, with the increase in number of constraints, the complexity of the solver increases proportionally. The complexity can be estimated by following the constraint set formulation and its properties in Section 3.1 of \cite{paper0}.

For DG discretizations, the same idea can be applied, except that straightforward mass conservation couples all elements. To retain the elementwise decoupling efficiency, we impose a stronger condition that the mass \emph{per element} remain unchanged. This allows us to use the same procedure as in \Cref{sec:fluxcons}, performing $E$ size-$n$ optimizations. In particular, we can simultaneously enforce both element boundary preservation and elementwise mass conservation with three equality constraints.

Note that these procedures are generalizable for arbitrary linear constraints. In particular, if we have $K$ linear constraints, then the dimension of the optimization problem can be reduced to an $(n-K)$-dimensional problem using the procedure described in Section \ref{sec:fluxcons}. This, in turn, means that this optimization is meaningful only if the original size of the problem is more than $K$. For example, for DG simulations we can only accommodate $n-1$ linear constraints per element if there are $n$ degrees of freedom per element.

%% file: Algorithm.tex
\section{Algorithm}
\label{sec:algo}

This section describes some algorithmic considerations of our approach. We emphasize that although in previous sections we have described details in different ways for CG versus DG discretizations, the fundamental tool, i.e., the algorithms in Section \ref{sec:summ}, are identical. The differences manifest only when we need to fit these types of discretrizations into the general framework of that section. This section will also highlight some algorithmic differences between the CG and DG discretizations that surface in the use of our optimization filter.

Note that, regardless of the discretization used, this filter preserves $L^2$ stability properties since it is norm-contractive. For example, if a CFL condition for $L^2$ stabilty of an unconstrained solver is used to determine a stable time-step value, then a filtered version of this solver will not require a change in time-step value to maintain this stability property.

\subsection{Algorithmic ingredients}

As a first step of the optimization procedure, we must transform any given basis $\{\phi\}_{j=0}^{N-1}$ that is assumed by the PDE solver to an orthonormal basis $\{\psi\}_{j=0}^{N-1}$. We use the standard procedure that requires application of the inverse of Cholesky decomposition of the mass matrix. In DG simulations, we can leverage elementwise decoupling as described in Section \ref{sssec:dg} so that the appropriate matrix algebra is performed locally on a local mass matrix associated with only a single element. However, when a continuous Galerkin formulation is used, the global mass matrix must be inverted, and this cost significantly slows down the optimization process. 

In CG, the structure of the global mass matrix $\bs{M}$ when using standard hat and bubble basis functions is ``nearly" diagonal, but the inverse of $\bs{M}$ (or its Cholesky factor) is a dense matrix. In our simulations, we precompute the required global matrices for orthonormalization of the basis elements and coordinates. This procedure is feasible for our problems in one spatial dimension. However, this process is no longer feasible if the global number of degrees of freedom $N$ becomes too large, as would happen with problems in two or three spatial dimensions.
%This is one way to reduce the cost to apply the filter to large systems which yield huge global matrices. However, a more elegant method to find $M$, which makes it computationally feasible to invert, is warranted here.

\subsubsection{DG specializations}

For DG discretizations, the global filtering operation is decoupled to act on individual elements. This element-by-element application 
%the elements of the domain can independently fix the coefficients of projection such that the constraints can be imposed on every element simultaneously. This element-by-element application of filter 
lends itself to parallel implementation. In addition, we need not call the filter on elements where the constraints are already satisfied. For example, if positivity is the constraint, then at every timestep we can flag elements where positivity may be violated and run the optimizer only over those elements. In our implementation, we use confederate linearization to check for the optima of function values on each element. This approach guarantees that all the optima are calculated accurately. We use the optima to determine whether there is a violation of positivity constraint on that element. Although a computational cost is associated with the eigen value solve on each element, the overall speed-up obtained by not having to calculate the signed distance function described in \eqref{eq:global-minimization} on each element certainly justifies the extra computation.

With this flagging scheme, we observe that the computational time spent on the filter is a fraction of the cost of the unconstrained solver (cf. our numerical results in Section \ref{sec:results}).

\subsection{Algorithm summary}

For completeness, we now summarize one procedure proposed in \cite{paper0} to perform the optimization in \Cref{sec:method}. Consider a stable implementation of \eqref{eq:contadr}. As discussed in \Cref{sec:method}, the optimization that enforces a constrained solution is performed every timestep.

Recalling some notation from Section \ref{sec:summ}, we assume $K$ families of linear constraints, and for a fixed $k \in \{1, \ldots, K\}$ each constraint is of the form \eqref{eq:constraints}. This results in $K$ feasible sets $\{C_k\}_{k=1}^K$, and $K$ signed distance functions $\{s_k\}_{k=1}^k$ each of the form \eqref{eq:sdist}. The full feasible set $C$ is the intersection of the $C_k$.

Computation of this signed distance function requires orthonormalization of a given basis for $V$,  which in turn requires the (inverse) Cholesky factor of the mass matrix, which can be precomputed before the simulation begins. A conversion from an input coordinate vector to coordinates $\tilde{\bs{v}}$ in an orthonormal basis is the first stage of the filtering procedure.

The second stage of the filter is an iterative procedure that attempts to project $\tilde{\bs{v}}$ onto the feasible set $C_k$.
At each iteration, we find the constraint index $k$ and the spatial point $x$ that minimizes the signed distance function, 
{
\begin{align}\label{eq:global-minimization}
 (x^\ast, k^\ast) \coloneqq \argmin_{x \in \Omega, k \in \{1, \ldots, K\}} s_k(x) %\mathrm{sdist}(\bs{\tilde{v}}, H_k(\tilde{v})),
\end{align}
%where $H_k$ defines the hyperplanes as shown in \Cref{fig:distandprojection}.
We then update $\bs{\tilde{v}}$ by projecting it onto the hyperplace corresponding to $(x^\ast, k^\ast)$:
\begin{align}\label{eq:greedy-update}
 \bs{\tilde{v}} \gets \bs{\tilde{v}} + \bs{h}(x^\ast, k^\ast) \min\left\{0, s_{k^\ast}(x^\ast)\right\}% \mathrm{sdist}(\bs{\tilde{v}}, H_{k^\ast}(\tilde{v}^\ast)) \right\}.
\end{align}
}
where $\bs{h}(x, k)$ is the normal vector corresponding to hyperplane $H_x$ of constraint family $k$ that points toward $C_k$. This vector is readily computable from the orthonormal basis, see \cite{paper0} for details. This procedure is repeated until $s_{k^\ast}(x^\ast)$ vanishes to within a numerical tolerance. Upon termination of the iterations, the output of the filtering procedure is the (updated) $\bs{\tilde{v}}$.

A brief summary of steps taken by a filtered PDE solver is presented in \cref{alg:greedy}.

\begin{algorithm}[H]
 \caption{Constrained PDE timestepping}
\label{alg:greedy}
\begin{algorithmic}[1]
\STATE{Input: Terminal time $T$, timestep size $\Delta t$, PDE solver spatial basis $\bs{\phi}$}
\STATE{Orthonormalize the basis function $\bs{\phi}$ to $\bs{\psi} \in \R^N$}
\STATE {Define nsteps $= \frac{T}{\Delta t}$}
\FOR{$i= 0,\cdots,$ nsteps }
\STATE{Solve PDE to obtain the coefficients $\bs{\tilde{v}}^{i} \in \R^N$}
\STATE{Input: constraints $(L_k, \ell_k)_{k=1}^K$}
\WHILE{True}
 \STATE{Compute $(x^\ast, k^\ast)$ via \eqref{eq:global-minimization}.}\label{alg:greedy:minimization}
 \STATE{If $s_{k^\ast}(x^\ast) \geq 0$, \textbf{break}}\label{alg:greedy:termination}

 \STATE\label{alg:greedy:update}{Update $\bs{\tilde{v}}^{i+1}(t)$ via \eqref{eq:greedy-update}.}
\ENDWHILE
\STATE {$\bs{{v}}^{i+1} = \bs{\tilde{v}}^{i+1}$}
%\RETURN {$\bs{\tilde{v}}^{i+1}$}
\ENDFOR
\end{algorithmic}
\end{algorithm}

%% file: Results.tex
\section{Numerical results}
\label{sec:results}

We numerically investigate the proposed procedure for preserving convex constraints in solutions to PDEs. We consider a bounded 1D spatial interval $\Omega \subset \R$. We will consider positivity constraints (i.e., $u(x) \geq 0$ over all $\Omega$), and will also investigate cell boundary value (``flux") preservation and mass conservation (for DG simulations). We are primarily interested in the effect that the filter has on convergence rates and in quantifying the computational efficiency of the procedure. All simulations perform the procedure summarized in Algorithm \ref{alg:greedy}. We use the following machine to report all the performance and error numbers in this section: 256 Intel(R) Xeon(R) CPU E7-4850 v4 @ 2.10GHz cores (HT) with 1024 GB of RAM running Redhat Enterprise 7.5 (Maipo).
%We present the application of the filter to the following types of PDEs in 1D and show empirical proof of convergence.

%\begin{itemize}	
% \item Discontinuous finite element experiment to solve 1D advection problem with varying constraints.
% \item Continuous finite element experiment to solve 1D diffusion reaction problem with positivity constraint.
%\end{itemize}

\subsection{DG}\label{sec:dgresults}

Consider the 1D advection equation,
\begin{align*}\label{eq:1dadv}
\frac{\partial {u}}{\partial t} + a \frac{\partial {u}}{\partial x} = 0,
\end{align*} 
where $a$ is a fixed constant. For $a = 1$, we investigate this problem for the exact solution ${u}(x,t) = 0.5 \sin(2\pi {x}- 2 \pi t-0.5 \pi) + 0.5$. We use a DG formulation with periodic boundary conditions over the domain $\Omega=[-1,1]$, which results in a system of ordinary differential equations prescribing the evolution of the Galerkin coefficient. We employ Runge-Kutta-2 to integrate in time and compute up to a final time $T=1$ using upwind flux calculation. %We consider unit advection velocity ($a=1$) for simplicity. 

%Because of periodic nature of the test function, we get the initial and final solution to overlap on $\Omega$. As such, an ideal final solution $u(x,T)$ would be the one that exactly matches the initial value $u(x,0)$ and respect the structure of $u(x,0)$. However, in practice, comparing $u(x,0) $ to $u(x,T)$, we find structural discrepancies. That is why it is an excellent choice of function to demonstrate the correction procedure by applying the filter.

Once fully discretized, the numerical solution can correspond to a function that is negative in some parts of the spatial domain, and hence we apply the filtering procedure in Algorithm \ref{alg:greedy} to enforce positivity on all $\Omega$. This results in two numerical solutions:
\begin{itemize}
 \item $\tilde{\bs{v}}$ is the solution resulting from a \emph{standard} DG solver, i.e., one that does \emph{not} employ our filter.
 \item $\bs{v}$ is the numerical solution resulting from application of the filter as specified in Algorithm \ref{alg:greedy}.
\end{itemize}
As discussed earlier, this filtered solution does not respect mass conservation, and it in general changes values on the boundary, which changes numerical fluxes. Therefore, we use Algorithm \ref{alg:greedy} to generate two more numerical solutions:
\begin{itemize}
 \item $\bs{v}_F$ is the positivity-constrained solution that adds in elementwise equality constraints to preserve inter-element boundary (flux) values.
 \item $\bs{v}_{I + F}$ is the positivity-constrained solution that includes both inter-element flux preservation and elementwise mass conservation.
\end{itemize}
%We apply the filter (Algorithm \ref{alg:greedy}) to enforce positivity, resulting in the solution that we label $\bs{v}$. 

\begin{figure}[htbp]
 \begin{center}
  \resizebox{\textwidth}{!}{
   \includegraphics[width=0.49\textwidth]{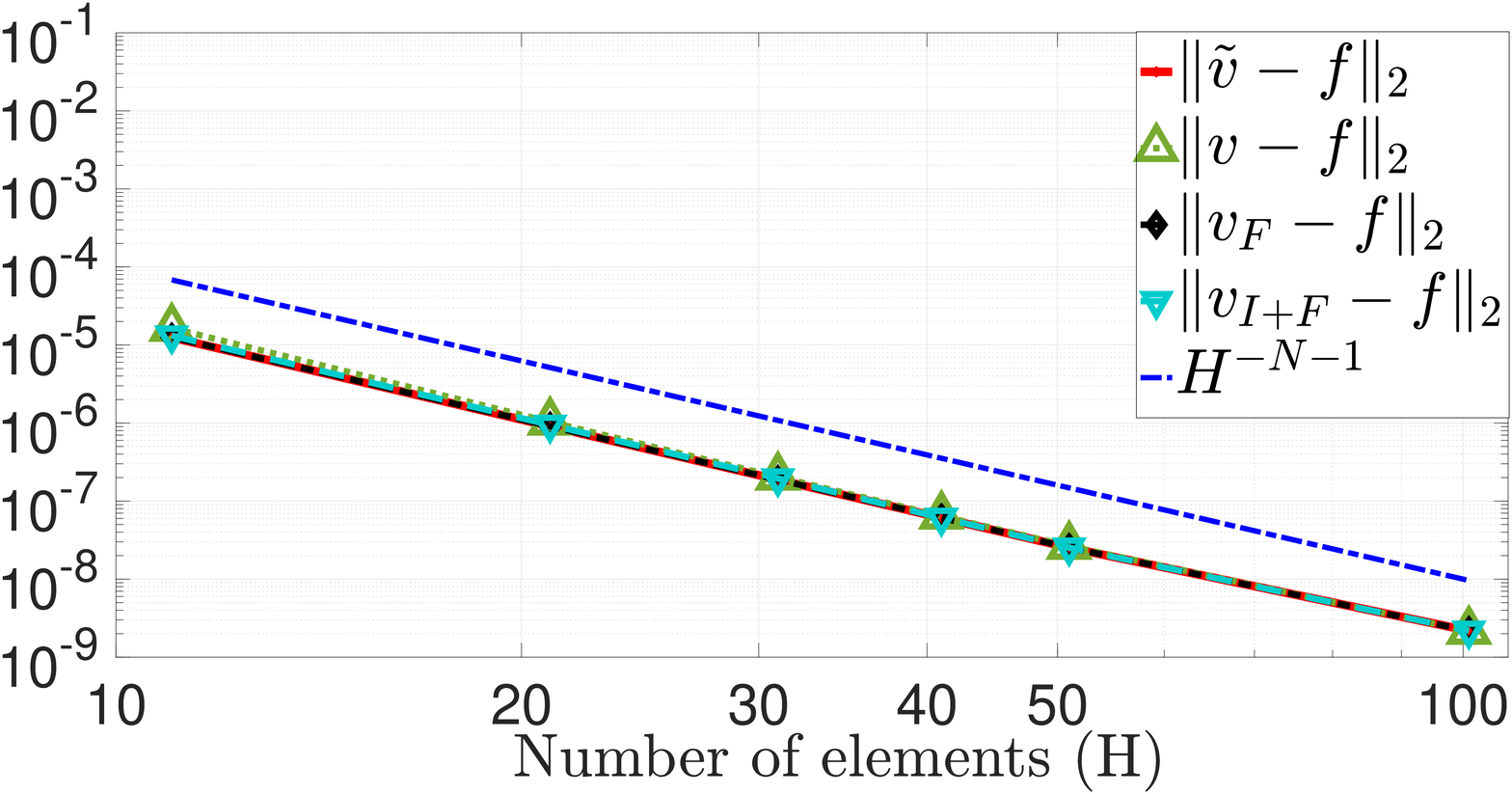}
   \includegraphics[width=0.49\textwidth]{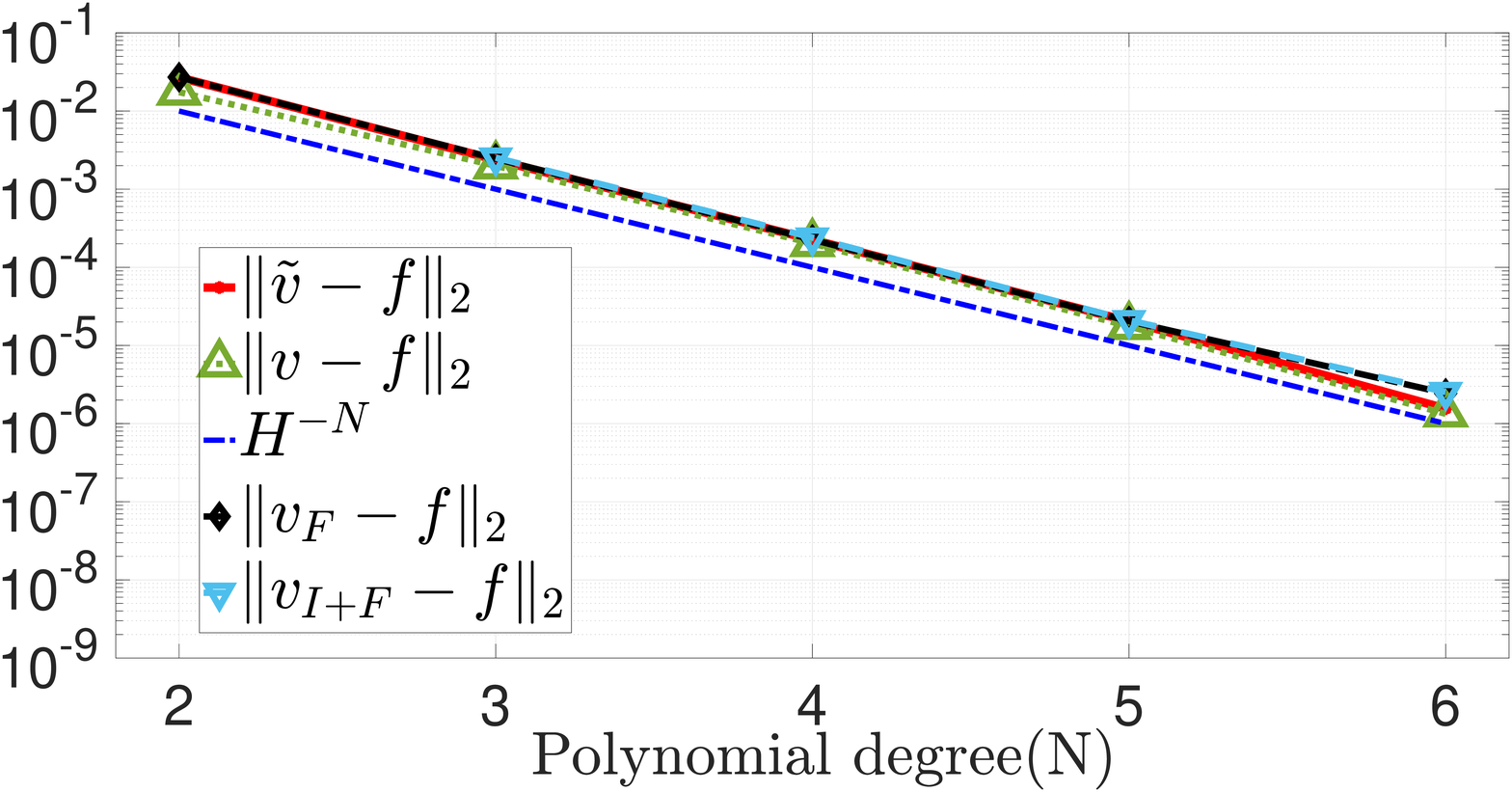}
  }
 \end{center}
 \caption{Convergence study for DG solution to the 1D advection PDE applied to ${u}(x,t) = 0.5 \sin(2\pi {x} - 2 \pi t-0.5 \pi) + 1$ for a time period of $T = 1$ second. $\bs{\tilde{v}}$ refers to the unfiltered solution, ${\bs{{v}}}$ refers to the positive solution, ${\bs{{v}}}_F$ refers to the positive solution with flux (boundary values) conserved for each element, and ${\bs{{v}}}_{I+F}$ refers to the positive solution with flux and mass conserved for each element . 
  Left: h-convergence using constant polynomial order $N = 3$, $\Delta t = 10^{-5}$. 
  Right: p-convergence at constant number of elements $H = 3$, $\Delta t = 10^{-4}$. } \label{fig:dgex1}
\end{figure}

%To solve by DG formulation, first we need to project the ${u}(x,t)$ onto the polynomial space using appropriate basis functions ${\phi}^N(x)$ such that ${u}(x,t) = \sum\limits_{i= 0}^{N-1} \tilde{v}_i^t\phi(x)$. Let us call $\{\tilde{v}_0, \cdots, \tilde{v}_{N-1}\} $ at time $t$ as $ \bs{\tilde{v}}^t$. Using $\bs{\tilde{v}}^t$ to evaluate ${u}(x,t) $, we encounter issues such as the evaluated solution becoming negative at certain points in the domain. To overcome this, we apply the filtering procedure on ${\bs{\tilde{v}}}^t$ to obtain $\bs{{v}}^t$ such that it evaluates to positive values on the entire continuous bounded domain $\Omega$. The side-effects of this as explained in \Cref{sec:fluxcons,sec:energycons} can be handled by including additional flux and energy conservation constraints to the existing positivity constraint set of the filter. 
\Cref{fig:dgex1} investigates both $h-$ and $p-$convergence of the numerical solution to \eqref{eq:1dadv} for all four numerical solutions. We observe that, regardless of the constraint that we impose, the convergence rates are unchanged, and even the error values are nearly identical for all numerical solutions. Thus, for this example, our proposed procedure can guarantee positivity (and flux/mass conservation as well, if desired) without a notable impact on the accuracy of the solver.

\begin{figure}[H]
 \centering
 \label{fig:pHat}
\includegraphics[width=0.6\textwidth]{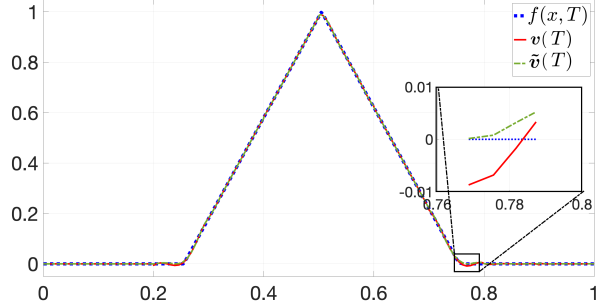}   \caption{Periodic function $f$ for convergence study of filtered 1D DG solution to \eqref{eq:1dadv}. Note that the $f(x,0)$ and $\bs{\tilde{v}}(T)$ overlap because of the periodic nature of $f$; however $\bs{\tilde{v}}(T)$ does not comply to the positivity structure of $f(x,0)$ as expected. After the application of filter, we obtain $\bs{{v}}(T)$, which changes $\bs{\tilde{v}}(T)$ to preserves the positive structure of $f$. }
\end{figure}

\begin{figure}[htbp]
 \begin{center}
  \resizebox{\textwidth}{!}{
   \includegraphics[width=0.49\textwidth]{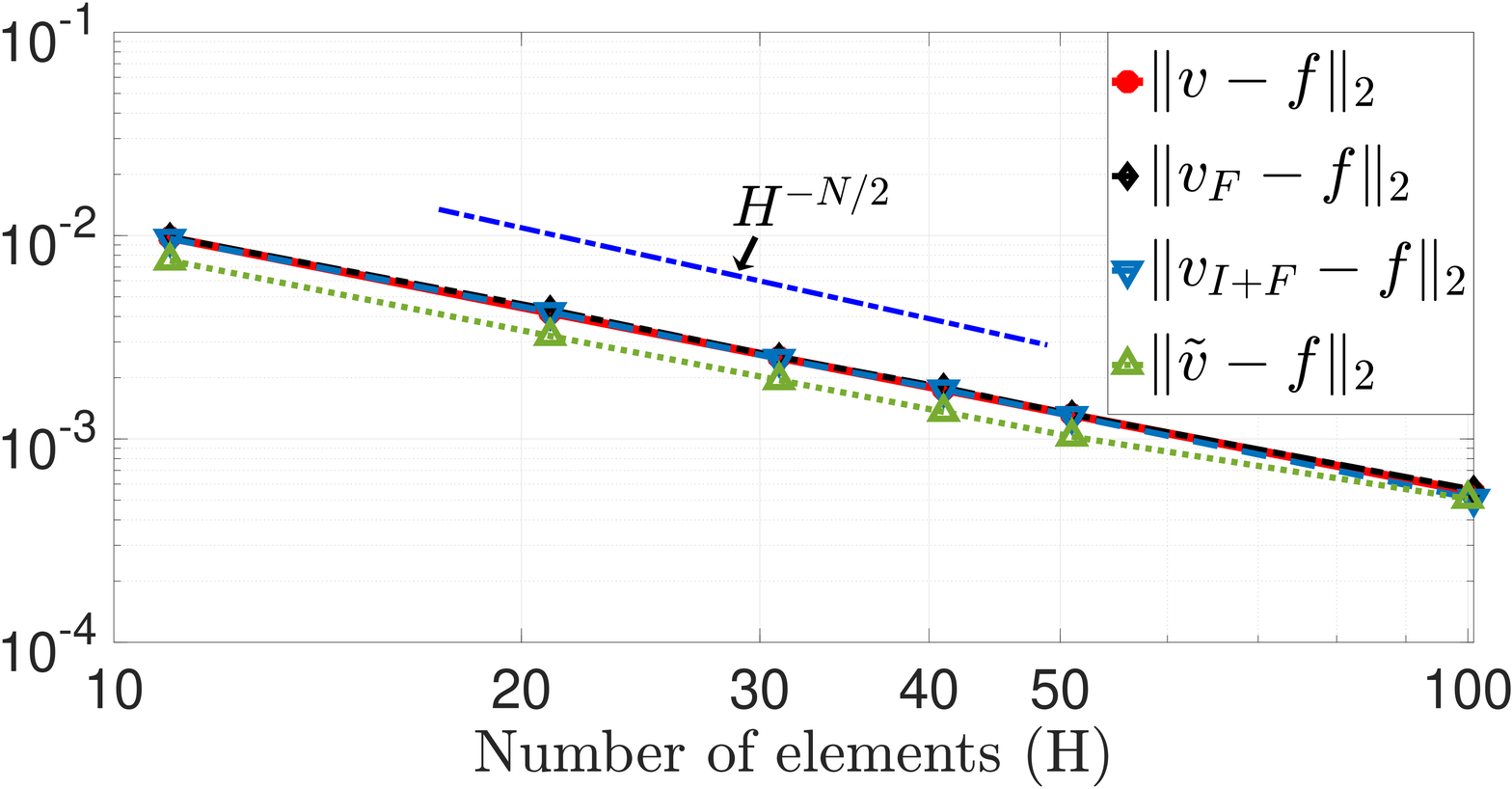}
   \includegraphics[width=0.49\textwidth]{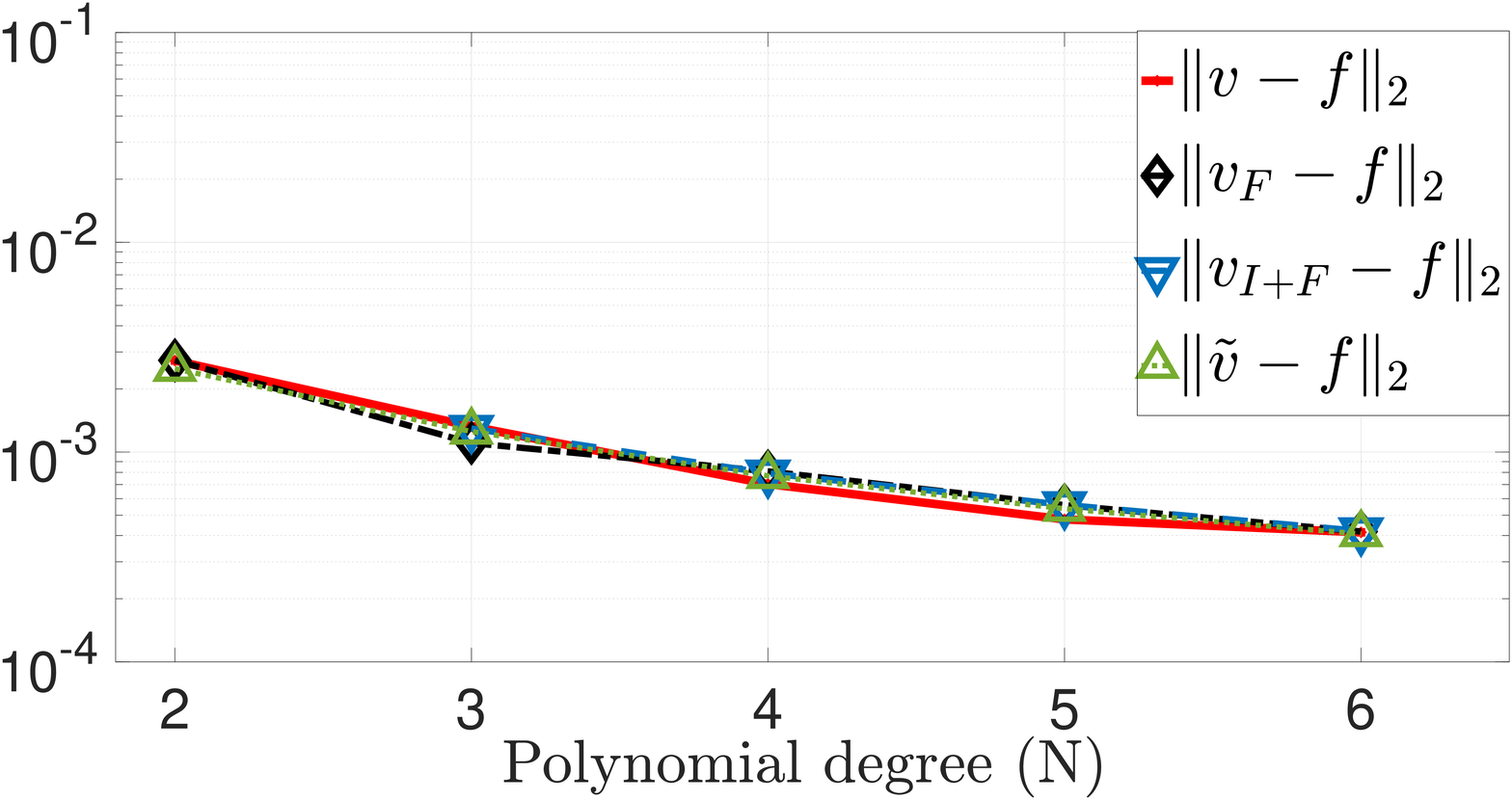}
  }
 \end{center}
 \caption{Convergence study for DG formulation of 1D advection PDE applied to \Cref{fig:pHat}, for a simulation run to T$=1$. Vector $\bs{\tilde{v}}$ represents the unfiltered solution coefficients, ${\bs{{v}}}$ refers to the positive solution coefficients at $T = 1$, ${\bs{{v}}}_F$ refers to the coefficients of the positive solution with flux (boundary values) conserved for each element, and ${\bs{{v}}}_{I+F}$ refers to the positive solution coefficients with flux and mass conserved for each element . 
  Left: h-convergence at constant polynomial degree $N= 3$ , $\Delta t = 10^{-5}$.
  Right: p-convergence at constant number of elements $H = 51$ , $\Delta t = 10^{-4}$. Here the choice of an odd number of elements for the $h-$convergence study ensures that the middle discontinuity of the function's derivative is located on an element boundary.} \label{fig:dgex2}
\end{figure}

We repeat the same convergence study for the same PDE, but instead for a triangular hat function with non-negative values as shown in \Cref{fig:pHat}; also shown in the figure is the unconstrained numerical solution at $T = 1$, which violates positivity. We run the same experiment as for the previous example and plot the $h-$ and $p-$convergence results in \Cref{fig:dgex2}.

From the results shown in \Cref{fig:dgex1,fig:dgex2} we observe that the convergence for filtered solution ${\bs{{v}}}$ and its variants remains largely unchanged/comparable to the unfiltered counterpart. Regarding the cost of the filtering procedure, there is a one-time orthonormalization cost of the basis function used as the initial setup, but since we can employ the filter over each element individually, this cost is negligible. The computational cost of filtering procedure per timestep depends on the time taken by the global minimum finding step of the filter, which we investigate next.

\begin{figure}[htbp]
 \begin{center}
  \resizebox{\textwidth}{!}{
   \includegraphics[width=0.49\textwidth]{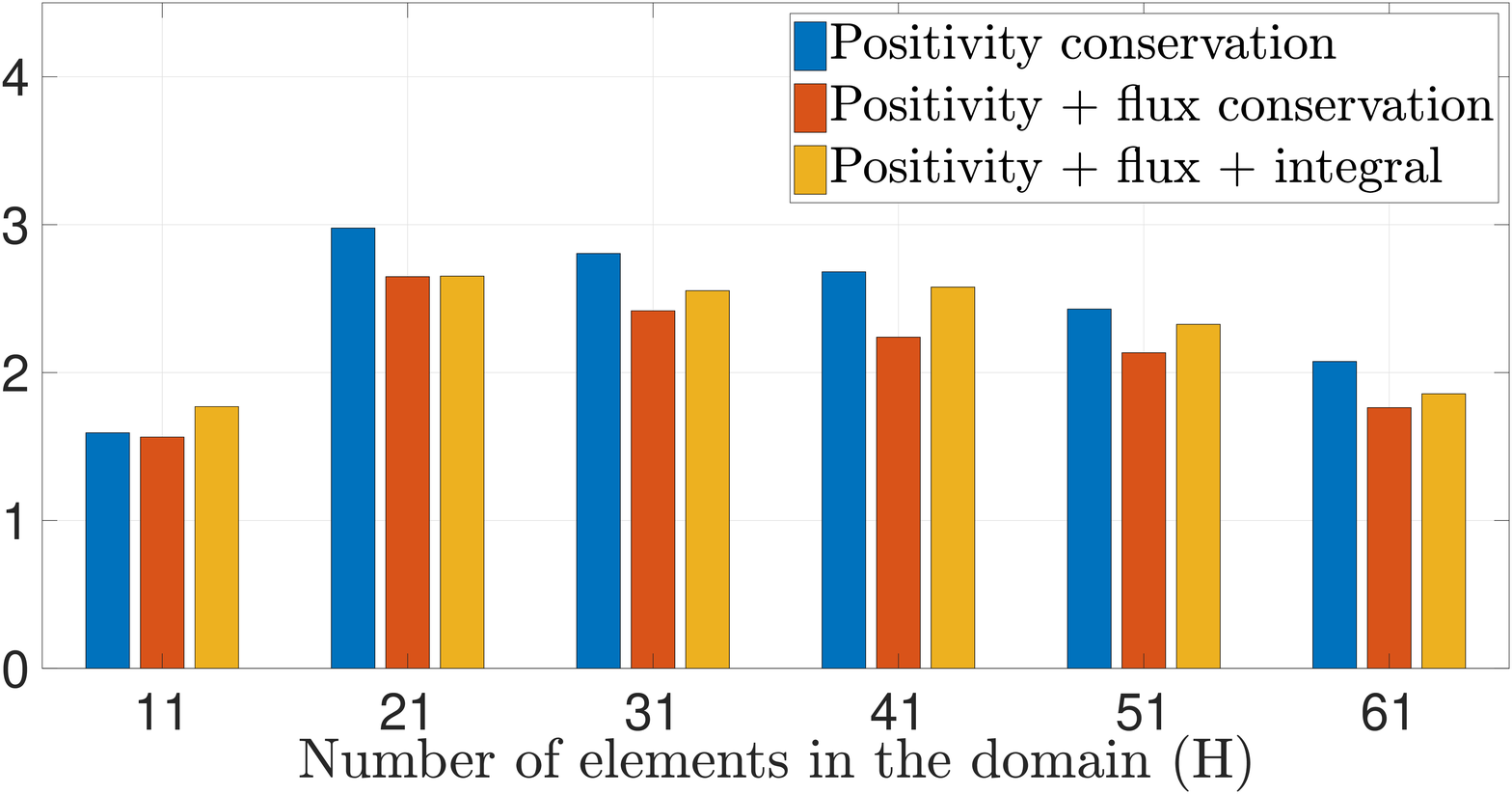}
   \includegraphics[width=0.49\textwidth]{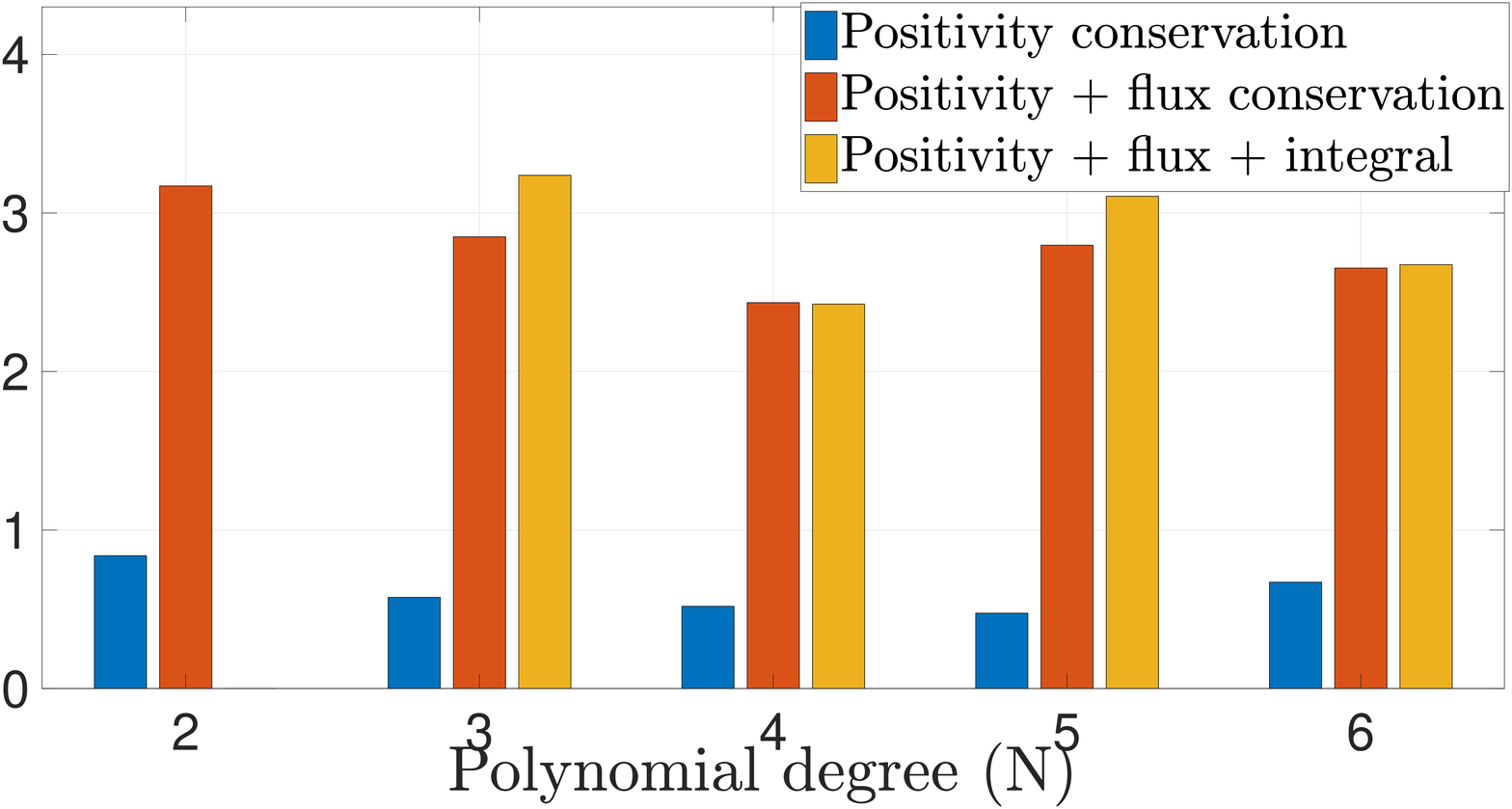}
  }
 \end{center}
 \caption{Average number of elements flagged per timestep during the experiments shown in \Cref{fig:dgex2}. Left: $h-$convergence experiment with $\Delta t = 10^{-5}$ and $N = 3$. Right: $p-$convergence experiment with $\Delta t = 10^{-4}$ and $H = 51$.} \label{fig:convhatavgflagged}
\end{figure}

For each simulation involving the filter, we compile the number of elements per timestep where the filter is employed. (Recall that we perform a cheaper check to certify positivity before calling the filter optimization.) For each of the three filtered solutions $\bs{v}$, $\bs{v}_F$, and $\bs{v}_{I+F}$, \Cref{fig:convhatavgflagged} plots the average number of filtered elements per timestep during the convergence experiments run in \Cref{fig:dgex2}. 

\begin{figure}[htbp]
 \begin{center}
  \resizebox{\textwidth}{!}{
   \includegraphics[width=0.49\textwidth]{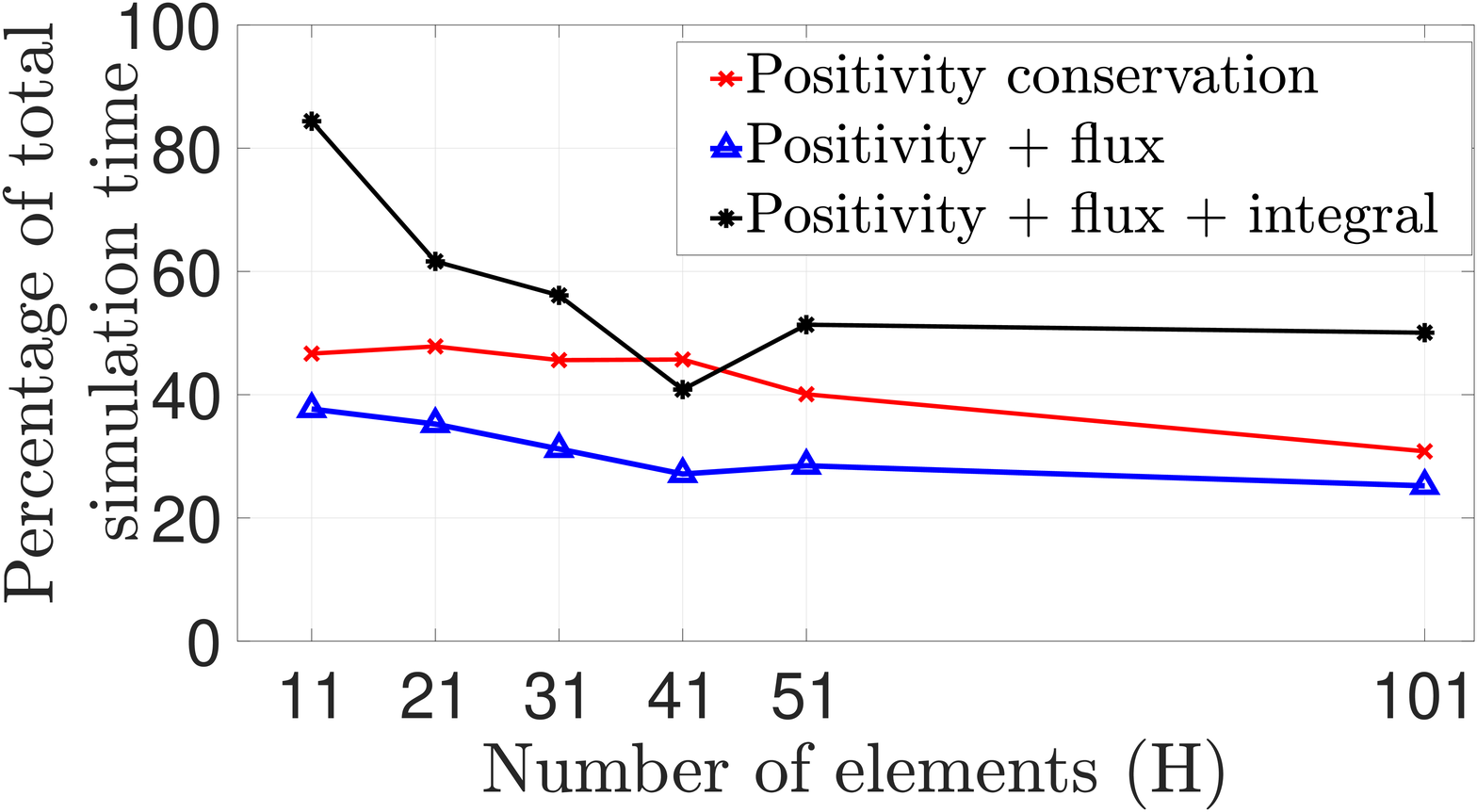}
   \includegraphics[width=0.49\textwidth]{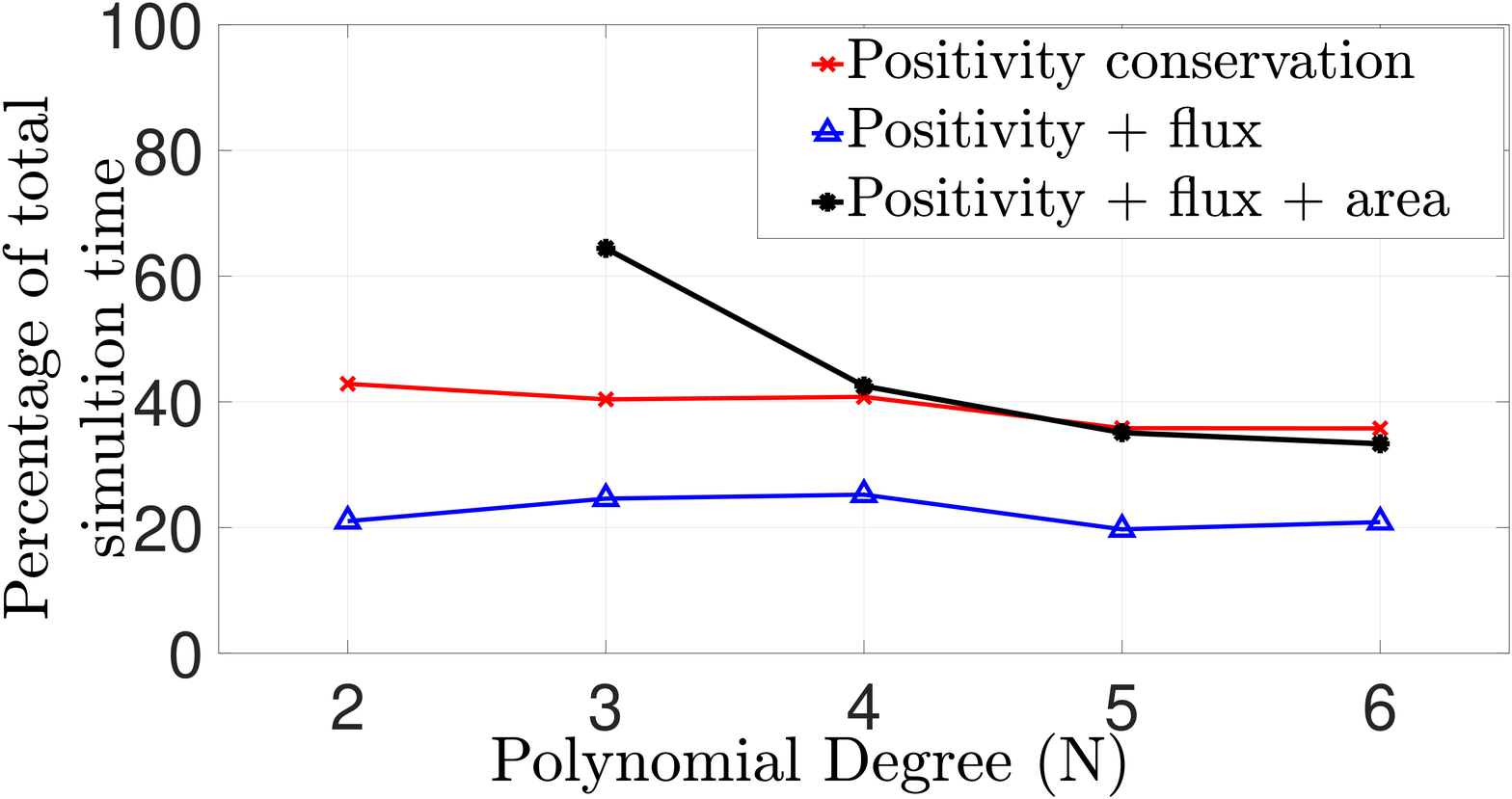}
  }
 \end{center}
 \caption{Percentage of total simulation clock time spent inside the filtering procedure during the convergence experiments shown in \Cref{fig:dgex2}. Left: $h-$convergence experiment with $\Delta t = 10^{-5}$ and $N = 3$. Right: $p-$convergence experiment with $\Delta t = 10^{-4}$ and $H = 51$.} \label{fig:perctime}
\end{figure}

As seen from \Cref{fig:convhatavgflagged}, the filter is not indiscriminately applied across all elements in the domain, and instead is called only for a few flagged elements per timestep where constraint violations are observed. This can substantially reduce the required computational overhead compared to applying the filter across every element. The percentage of total simulation time that is spent inside the filtering procedure for experiment \Cref{fig:dgex2} is shown in \Cref{fig:perctime}.

For some experiment configurations, the filtering cost can be more than half the simulation time (e.g., $H=11$ for the black line in \Cref{fig:perctime}, left). However, in many cases, the filtering procedure can require less than $20\%$ of the total computational effort.
It is interesting to note that for the $p$-convergence experiment in \Cref{fig:perctime}, the percentage of filtering time is \emph{higher} for some filters with \emph{fewer} constraints. For example, in \Cref{fig:perctime}, right, less simulation time is required to impose positivity and flux preservation compared to simply imposing positivity. One possible explanation for this phenomenon is that filters with a larger number of equality constraints require optimization in a lower dimensional space. Thus, the optimization problem can be less expensive to solve in this case. However, with additional constraints (e.g., the black line in \Cref{fig:perctime}, right), even though the dimensionality of the optimizer decreases, the number of iteration required by the filter increases. Therefore, a suitable balance must be struck from application to application.

\subsection{CG}

We now present the results on the continuous Galerkin implementation of a diffusion-reaction equation, 
\begin{align*}\label{eq:cgdr}
\frac{\partial}{\partial t}{u}({x},{t}) = \gamma \frac{\partial^2}{\partial x^2} {u}({x},{t}) + r(u({x},{t})).
\end{align*}
We consider the problem on a bounded domain $\Omega \subset \R$ and $\Omega = [-1,1]$. Consider the following advecting smoothed Heaviside function as a solution to this problem:
\begin{align*}\label{eq:exactsolcg}
u(x,t) = e^{-\gamma t} \Big(tanh(\epsilon (x+0.4) - ct)+1\Big),
\end{align*}
where $\epsilon$ is a shape-parameter, and $c$ refers to the speed at which the exact solution moves in space over time. The nonlinear reaction term in this experiment is set as $r({u(x,t)}) = \mu {u(x,t)} (1-{u}^2(x,t))$, where $\mu$ is a constant. Our Galerkin discretization of this problem integrates $r$ exactly, which is possible since $r$ depends polynomially on $u$.
%Since $r(x,t)$ is a quadratic in $u(x,t)$ with degree $N=3 N$, we use galerkin projection using basis $\{\phi_0, \cdots, \phi_{3 N}\}$ to obtain $\bs{r} =\{\hat{r}_0,\cdots,\hat{r}_{3 N - 1}\}$. 
A diffusion-reaction PDE with a stiff quadratic reaction term is interesting as it represents a source term with the potential to blow up and only becomes active at the discontinuity \cite{leveque1990study}. Whereas numerical schemes may converge, they may converge with approximants that violate positivity.
%This is prone to make the numerical solution lose the desirable structural properties.

Using the discrete form of the continuous problem \eqref{eq:contadr} as described in \Cref{sec:method}, we choose the IMEX CNAB-2 scheme, which deals with the advection term using Crank-Nicolson, and uses a second-order Adams-Bashforth scheme for the reaction term.
The \emph{unconstrained} CG discrete time update then is
\begin{align*}
\bs{\tilde{v}}^{n+1} = \bs{B}^{-1}\Big(\bs{A}\bs{\tilde{v}}^{n} + \bs{{r}}^{n}\Big) 
\end{align*}
where $\bs{A}$ and $\bs{B}$ are given by
\begin{align*}
\bs{A} = \bs{M} - (0.5 \gamma \Delta t \bs{L}) \\
\bs{B} = \bs{M} + (0.5 \gamma \Delta t \bs{L}),
\end{align*}
with $\bs{M}$ and $\bs{L}$ the mass and Laplacian matrices, respectively, and $\bs{{r}}^n$ the projection coefficients of the reaction function $r(\bs{\tilde{v}}^n)$. 

\begin{figure}[H]
 \centering
 \label{fig:cGfunc1}\includegraphics[width=0.6\textwidth]{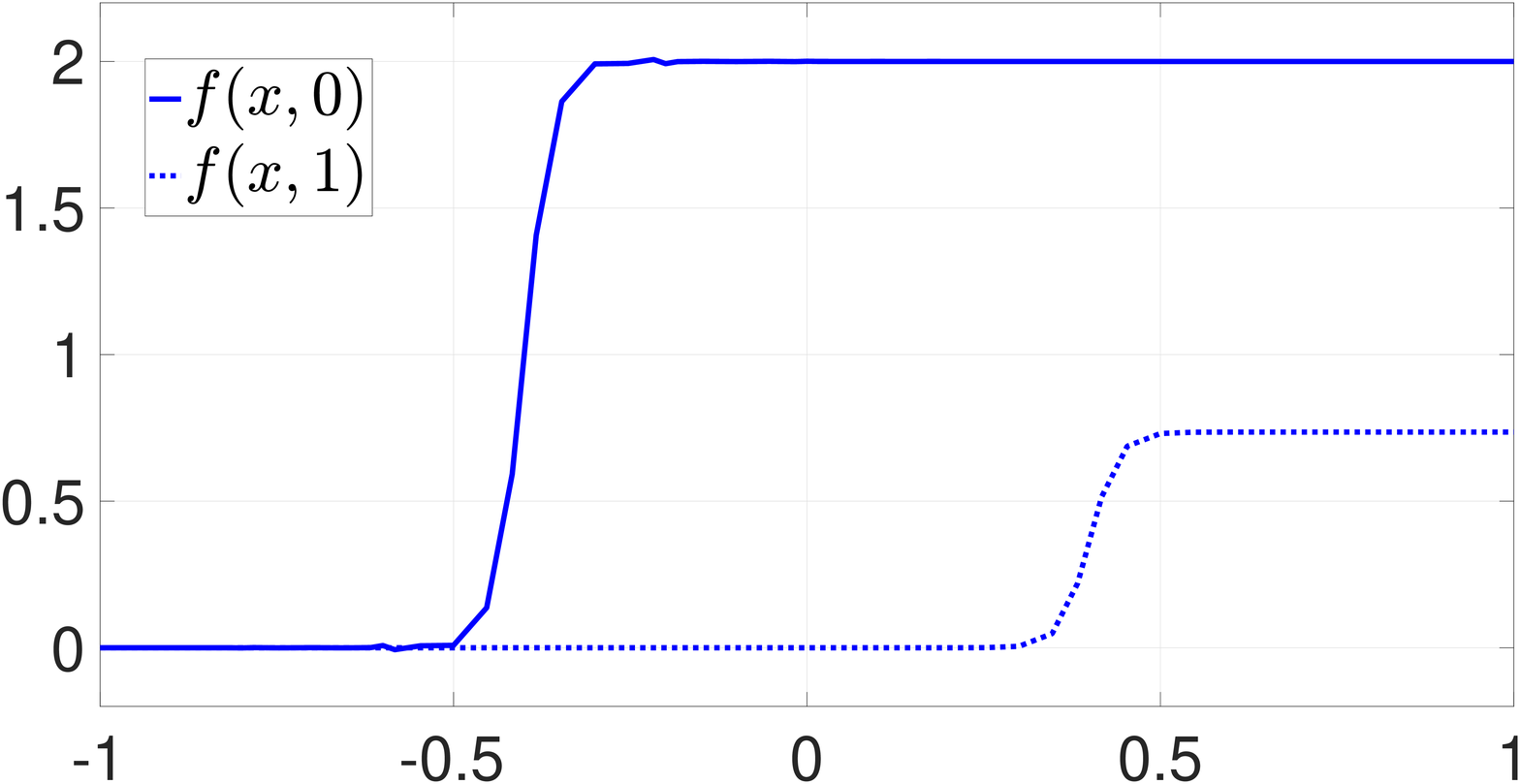}
 \caption{Initial and final state of the function $u$ for $\epsilon = 25$, $c=20$, and $\gamma = 1$, $\mu = 1$.}
\end{figure}

\begin{figure}[htbp]
 \begin{center}
  \resizebox{\textwidth}{!}{
   \includegraphics[width=0.49\textwidth]{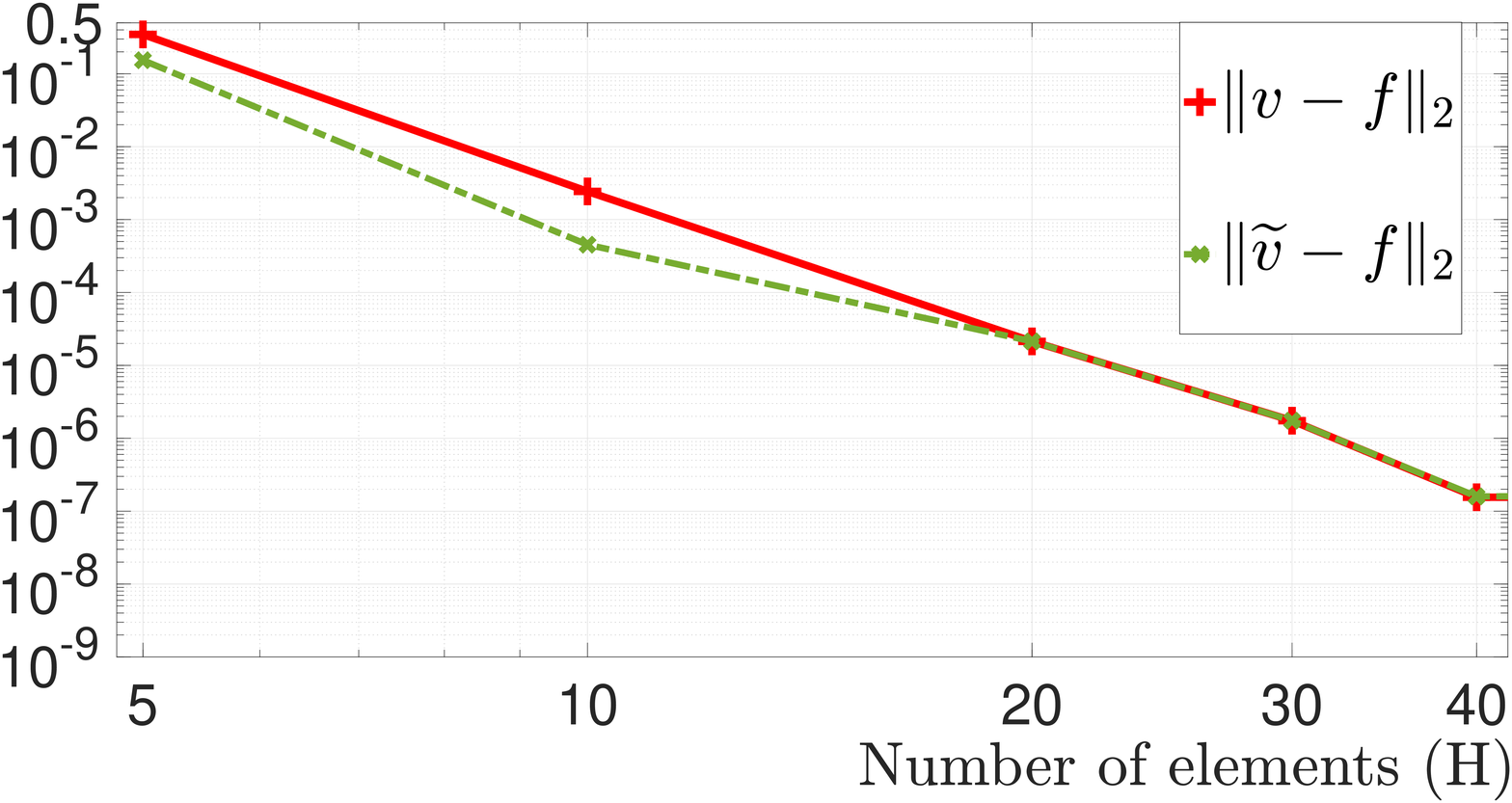}
   \includegraphics[width=0.49\textwidth]{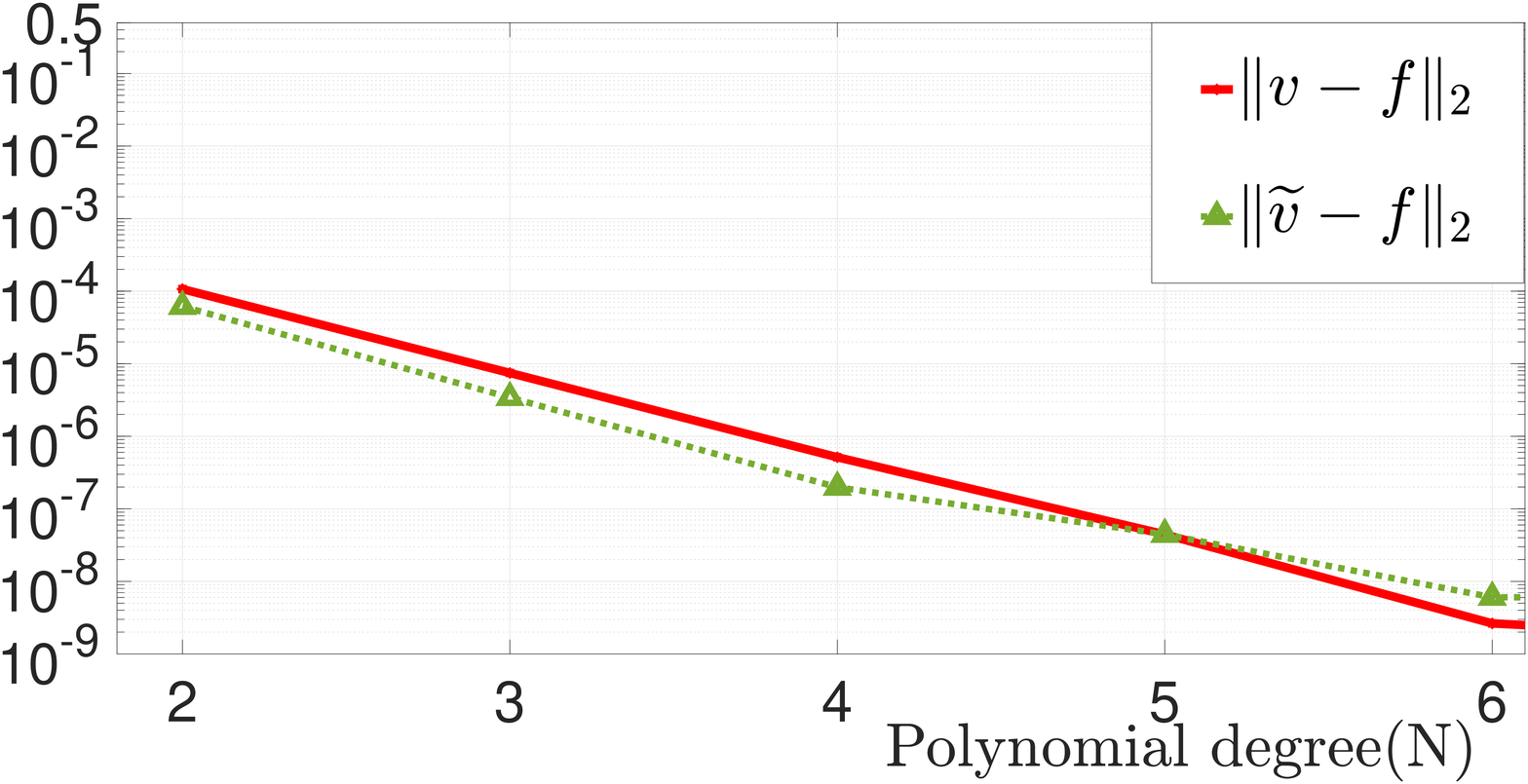}
  }
 \end{center}
 \caption{
 Convergence study for 1D CG diffusion-reaction PDE applied to the function $f$ shown in \Cref{fig:cGfunc1} with constant timestep $\Delta t = 10^{-4}$ until final time $t = 1$.
  Left: h-convergence at constant polynomial degrees $N= 7$.
  Right: p-convergence at constant number of elements $H = 100$.} \label{fig:cgexpt1}
\end{figure}

%For this experiment, we choose a second-order IMEX scheme CNAB-2: Crank Nicolson for advection and second-order Adams Bashforth for reaction. 
\Cref{fig:cGfunc1} shows the initial and final states of $u(x,t)$ for a simulation up to $T=1$. As described earlier, the CG formulation uses a nonorthonormal (``hat" and ``bubble" function) basis, meaning that a spatially global mass matrix inversion must be performed for the filtered versions of this experiment. Advances that make this procedure efficient are currently being investigated, but \Cref{fig:cgexpt1} demonstrates that, like the DG formulation, error and convergence rates are largely unaffected for this problem by inclusion of the filter. If the constraints are to be satisfied locally in the CG formulation, the order of accuracy we observe in this example may not hold. This problem is widely known in the context of flux limiters, and many fixes have been proposed, e.g., \cite{l_1}.

This example also demonstrates the flexibility of our approach: Although CG and DG discretizations can be very different mathematically and computationally, Algorithm \ref{alg:greedy} is the generic template for applying our procedure, independent of most PDE discretization details. Note that we consider global constraint satisfaction in the current example. In the case in which an additional advection term is present in the PDE, it conserving the structure of the solution locally is often desirable. Such conservation can be achieved by the current method after reformulating the optimization problem such that E constraints are in the constraint set, where E is the number of elements in the domain. Each constraint represents the property to be conserved in each element. 

%% file: Conclusions.tex
\section{Conclusions}
\label{sec:conclusions}

We have proposed a filtering approach applied to standard finite element time-stepping discretizations of PDE, using both continuous and discontinuous Galerkin formulations in one spatial dimension. The goal of the filter is to enforce one or more constraints from a general class of convex inequalities that can apply over the entire spatial domain. Thus, we pose the filtering problem as a nonlinear optimization problem, which we propose to perform after every time step. %We consider the problem of structure preservation of polynomial projection on each element during the timestepping procedure. 
We focus mainly on enforcing positivity, but enforcing a maximum principle or monotonicity also falls into our framework. Linear equality constraints, such as total mass preservation or elementwise boundary value preservation, are possible in our framework. Our numerical results show that the filtering (optimization) does increase
in the computational cost of PDE solvers, but it often requires less than 50\% of
the total simulation time.

%% file: Acknowledgments.tex
\section*{Acknowledgments}

V. Zala and  R.M. Kirby acknowledge that their part of this research was sponsored by ARL under cooperative agreement number W911NF-12-2-0023. The views and conclusions contained in this document are those of the authors and should not be interpreted as representing the official policies, either expressed or implied, of ARL or the U.S. Government. The U.S. Government is authorized to reproduce and distribute reprints for Government purposes notwithstanding any copyright notation herein.  A. Narayan was partially supported by NSF DMS-1848508. This material is based upon work supported by both the National Science Foundation under Grant No. DMS-1439786  and the Simons Foundation Institute Grant Award ID 507536 while A. Narayan was in residence at the Institute for Computational and Experimental Research in Mathematics in Providence, RI, during the Spring 2020 semester.